\documentclass[11pt,a4paper]{drna}

\usepackage{longtable}

\usepackage{mathtools}
    \mathtoolsset{centercolon}
\usepackage[ruled,noend]{algorithm2e}
\SetKwInput{KwResultT}{Result upon Termination}
\usepackage[inline]{enumitem}
\usepackage{fancyvrb}
\usepackage{multicol}
\usepackage{placeins}
\usepackage{tikz}
    \usetikzlibrary{arrows, positioning}
\usepackage{upquote}
\usepackage{xtab}

\newcommand{\set}[1]{\left\lbrace #1 \right\rbrace}
\let\set\set
\newcommand{\XX}[1]{{\color{red}\footnotesize\framebox{\textbf{XX}}\emph{\textbf{#1}}}}
\newcommand{\vardot}{\mathord{\,\cdot\,}}
\newcommand{\JSR}{\operatorname{JSR}}
\newcommand{\len}{\operatorname{len}}
\newcommand{\co}{\operatorname{co}}
\newcommand{\coast}{\operatorname{co}_{\operatorname{\ast}}}
\newcommand{\cocone}{\operatorname{co}_{\operatorname{+}}}
\newcommand{\cosym}{\operatorname{co}_{\operatorname{s}}}
\newcommand{\coell}{\operatorname{co}_{\operatorname{e}}}
\newcommand{\closure}{\operatorname{cl}}

\newcommand\Tstrut{\rule{0pt}{2.6ex}}         
\newcommand\Bstrut{\rule[-0.9ex]{0pt}{0pt}}   

\xentrystretch{0.06}  

\renewcommand{\XX}[1]{\relax}  

\begin{document}

\title{The finiteness conjecture for $3\times 3$ binary matrices}

\author{Mejstrik}{Thomas}{a}

\affiliation{a}{University of Vienna, Austria, {e-mail: \texttt{\small  thomas.mejstrik@gmx.at}};
The author is sponsored by the Austrian Science Foundation (FWF) grant P~33352.}

\volume{7}
\firstpage{1}

\maketitle


\begin{abstract}
The invariant polytope algorithm was a breakthrough in the joint spectral radius computation,
allowing to find the exact value of the joint spectral radius for most matrix families~\cite{GP2013,GP2016}.
This algorithm found 
many applications in problems of functional analysis, approximation theory, 
combinatorics, etc.. 

In this paper we propose a modification of the invariant polytope algorithm 
enlarging the class of problems to which it is applicable.
Precisely, we introduce mixed numeric and symbolic computations.
A further minor modification of augmenting the input set with additional matrices
speeds up the algorithm in certain cases.

With this modifications we are able to automatically prove the finiteness conjecture for
all pairs of binary $3\times 3$ matrices and sign $2\times 2$ matrices.
\end{abstract}


\section{Introduction}
In this paper we are concerned with the maximal asymptotic growth rate of products of matrices,
the so called \emph{joint spectral radius}.
It has been defined in 1960~\cite{RS1960} and since found applications 
in many seemingly unconnected areas of mathematics and engineering,
e.g.\  for
computing the regularity of wavelets and of subdivision schemes~\cite{DL1992}, 
the capacity of codes~\cite{MOS01},
the stability of linear switched systems~\cite{Gur95}.

\begin{definition}
\label{def_jsr}
\XX{checked}
Given a finite set of matrices $\mathcal{A}\subseteq\RR^{s\times s}$. 
The \emph{joint spectral radius} (JSR) of $\mathcal{A}$ is defined as
\begin{equation}\label{equ_jsr}
\JSR(\mathcal{A}):=
\lim_{n\rightarrow \infty}
\ \max_{A_j\in\mathcal{A}}\ %
\norm{A_{j_n} \cdots A_{j_1}}^{1/n},
\end{equation}
where $\norm{\vardot}$ is any sub-multiplicative matrix norm.
\end{definition}

\noindent
An open question in the joint spectral radius theory is the so called \emph{finiteness conjecture}~\cite{LW1995}: 
\begin{center}
\itshape
Given a matrix set, does there exist a finite product whose\\
powers' spectral radii attain the growth rate equal to its joint spectral radius?
\end{center}
The finiteness conjecture has been proven false, 
in the sense that such a finite product does not always exist;
although this case seems to be exceptional~\cite{BM2002,BTV2003,Koz2005,JB2008,HMST2011}.
In this paper we proof the finiteness conjecture for pairs of binary matrices of dimension~3.

\subsection{Overview and main results}
In Section~\ref{sec_ipa} we present the \emph{invariant polytope algorithm (ipa)}
for computing the $\JSR$ of a finite set of square matrices.
In Section~\ref{sec_symbolic} we discuss how mixed numeric-symbolic computations 
can be used to widen the classes of matrices where the ipa is applicable.
In Section~\ref{sec_limitmatrix} we show how we can augment our input set of matrices
to obtain faster termination properties.

Finally, in Section~\ref{sec_fc} we discuss how the ipa, 
together with the discussed modifications,
can automatically proof the finiteness conjecture for pairs of binary matrices 
of dimension $2$ and $3$,
as well for pairs of sign matrices of dimension $2$.

\subsection{Notation}
The set of all integers is denoted by~$\NN$, 
integers including zero by~$\NN_0$, 
reals by~$\RR$, 
non-negative reals by~$\RR_+$,
complex numbers by~$\CC$.
Given $X\subseteq\CC^s$, where $s\in\NN$ is the dimension,
we denote the closure of $X$ by $\closure(X)$ and its interior by $X^\circ$.
Products of sets are understood element wise,
e.g.\ $A\cdot B=\set{a \cdot b:a\in A,\ b\in B}$.
Comparisons of matrices are understood element wise.
For a matrix $A$ we denote by $A^T$ its transpose 
and for a square matrix by $\rho(A)$ its spectral radius.

We will make use of various convex hulls of sets throughout the paper.
\begin{definition}\label{def_hull}
\begin{itemize}[leftmargin=*]
    \item For $V\subseteq\RR^s$, we define its \emph{convex hull} $\co V$ as the
    intersection of all convex sets containing $V$.
    
    
    \item  
    For $V\subseteq\RR^s_+$, we define the \emph{cone hull} of $V$ 
    (in the first orthant) 
    by
    \begin{equation}\label{equ_coast_cominus}
    \cocone V=
    \set{x\in \RR^s_+ : x=y-z,\ y\in\co(V),\ z\in\RR^s_+}\subseteq\RR^s_+.
    \end{equation}
    
    \item For $V\subseteq\RR^s$, we define the \emph{symmetric convex hull} of $V$  by
    \begin{equation}\label{equ_coast_cos}
    \cosym V = 
    \co\set{ V,\, -V }\subseteq\RR^s.
    \end{equation}    
    
    \item For $v=a+ib\in\CC^s$ we define its corresponding ellipse $E(v)=E(a,b)\subseteq\RR^s$
    as the two dimensional subset $\set{a \cos t + b \sin t:t\in\RR}\subseteq\RR^s$.
    For $V\subseteq\CC^s$, we define the \emph{elliptic convex hull}
    of $V$ by
    \begin{equation}\label{equ_coast_absco}
    \coell V = \co \set{ E(v) : v\in V}\subseteq\RR^s.
    \end{equation}

    \item     
    For simplicity, we denote with $\coast V$ any of the convex hulls 
    $\cocone$, $\cosym$, $\coell$,
    depending on the context.
\end{itemize}
\end{definition}

We will use the aforementioned convex hulls to define norms via their unit ball.
\begin{definition}
Let $P\in\RR^s$ be a compact, convex set with non-empty interior,
and such that
$rP\subseteq P$ for all $\abs{r}\leq1$.
We define the \emph{Minkowski norm}
$\norm{\vardot}_P:\RR^s\rightarrow\RR$ by
\begin{equation}
\label{equ_minknorm}
\norm{x}_P = \min \set{r>0:x\in r P}.
\end{equation}
\end{definition}

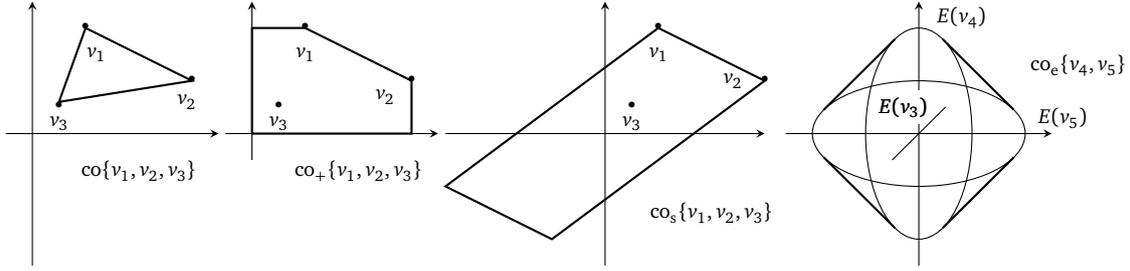
\begin{figure}
\centering
{\small
\begin{tikzpicture}[scale=0.70]
    \draw [-stealth](-.5,0) -- (3.5,0);
    \draw [-stealth](0,-2.5) -- (0,2.5);
    
    \node at (1.2,1.5) {$v_1$};
    \node at (1,2) {\textbullet};
    \node at (2.9,0.6) {$v_2$};
    \node at (3,1) {\textbullet};
    \node at (0.5,0.2) {$v_3$};
    \node at (0.5,0.5) {\textbullet};
    \node at (2,-0.7) {$\co \{v_1,v_2,v_3\}$};
    \draw [thick] (1,2) -- (3,1) -- (0.5,.6) -- (1,2);
\end{tikzpicture}
\begin{tikzpicture}[scale=0.70]
    \draw [-stealth](-.5,0) -- (3.5,0);
    \draw [-stealth](0,-0.5) -- (0,2.5);
    \draw [color=white](0,-2.5) -- (0,-0.5);
    
    \node at (1,1.5) {$v_1$};
    \node at (1,2) {\textbullet};
    \node at (2.5,0.7) {$v_2$};
    \node at (3,1) {\textbullet};
    \node at (0.5,0.2) {$v_3$};
    \node at (0.5,0.5) {\textbullet};
    \node at (2,-0.7) {$\cocone \{v_1,v_2,v_3\}$};
    \draw [thick] (0,0) -- (0,2) -- (1,2) -- (3,1) -- (3,0) -- (0,0);
\end{tikzpicture}
\begin{tikzpicture}[scale=0.70]
    \draw [-stealth](-3,0) -- (3,0);
    \draw [-stealth](0,-2.5) -- (0,2.5);
    
    \node at (1,1.5) {$v_1$};
    \node at (1,2) {\textbullet};
    \node at (2.4,1) {$v_2$};
    \node at (3,1) {\textbullet};
    \node at (0.5,0.2) {$v_3$};
    \node at (0.5,0.5) {\textbullet};
    \node at (2,-1.5) {$\cosym \{v_1,v_2,v_3\}$};
    \draw [thick] (1,2) -- (3,1) -- (-1,-2) -- (-3,-1) -- (1,2);
\end{tikzpicture}
\begin{tikzpicture}[scale=0.70]
    \draw [-stealth](-2.5,0) -- (2.5,0);
    \draw [-stealth](0,-2.5) -- (0,2.5);

    \node[fill=white] at (-.3,0.5) {$E(v_3)$};
    \node at (-.3,0.5) {$E(v_3)$};
    \draw    (-.5,-.5) -- (.5,.5);
    \draw    (0,0) ellipse (2 and 1);
    \node at (0.8,2.2) {$E(v_4)$};
    \draw    (0,0) ellipse (1 and 2);
    \node at (2.7,0.3) {$E(v_5)$};
    \node at (3,1.3) {$\coell \{v_4,v_5\}$};
    
    \draw  [thick] ( 1.788854381999832, 0.447213595499958) -- ( 0.447213595499958, 1.788854381999832);
    \draw [thick] (-1.788854381999832, 0.447213595499958) -- (-0.447213595499958, 1.788854381999832);
    \draw [thick] ( 1.788854381999832,-0.447213595499958) -- ( 0.447213595499958,-1.788854381999832);
    \draw [thick] (-1.788854381999832,-0.447213595499958) -- (-0.447213595499958,-1.788854381999832);
\end{tikzpicture}
}
\caption[Various convex hulls]%
{
Various convex hulls. 
$v_1=\begin{bsmallmatrix}1\\2\end{bsmallmatrix}$,
$v_2=\begin{bsmallmatrix}2\\1\end{bsmallmatrix}$,
$v_3=\begin{bsmallmatrix}0.5\\0.5\end{bsmallmatrix}$,
$v_4=\begin{bsmallmatrix}2\\i\end{bsmallmatrix}$,
$v_5=\begin{bsmallmatrix}i\\2\end{bsmallmatrix}$.
}
\label{fig_hull}
\vspace{-5mm}
\end{figure}

\section{The invariant polytope algorithm}
\label{sec_ipa}

The \emph{invariant polytope algorithm (ipa)} for the computation of the $\JSR$
makes use of the inequality~\cite{DL1992},
\begin{equation}\label{equ_bound_jsr}
\max_{A_{j}\in\mathcal{A}}\rho\left( A_{j_k}\cdots A_{j_1}\right)^{1/k}
\leq
\JSR(\mathcal{A})
\leq
\max_{A_j\in\mathcal{A}}
\big\| 
A_{j_k}\cdots A_{j_1}
\big\|^{1/k},
\end{equation}
which holds for any $k\in\NN$ and any sub-multiplicative norm $\norm{\vardot}$.
Before we can describe how the ipa works, we need a few further definitions.

\begin{definition}
\textbullet\ %
For a product $A_{j_k}\cdots A_{j_1}$ we say the number 
$\rho\left( A_{j_k}\cdots A_{j_1}\right)^{1/k}$ is its \emph{averaged spectral radius}. 

If there exists a product $\Pi=A_{j_n}\cdots A_{j_1}$, $A_j\in\mathcal{A}$, such that $\rho(\Pi)^{1/n}=\JSR(\mathcal{A})$, 
i.e.\ its averaged spectral radius equals the joint spectral radius,
we call the product a 
\emph{spectral maximizing product (s.m.p.)}.

\textbullet\ %
Given a matrix $\Pi\in\RR^{s\times s}$,
we call the eigenvalues largest in modulus the \emph{leading eigenvalues} 
and the corresponding eigenvectors the \emph{leading eigenvectors}.
If there exists only one largest eigenvalue in modulus (counted with algebraic multiplicity),
we say the leading eigenvalue is~\emph{simple}.

\textbullet\ %
Given a bounded set of matrices $\mathcal{A}\subseteq\RR^{s\times s}$ 
with $\lambda=\JSR(\mathcal{A})>0$
and let $\tilde{\mathcal{A}}$ be the set of normalized matrices
$\tilde{\mathcal{A}} = \{ A_j/\lambda : A_j\in\mathcal{A} \}$.
$\mathcal{A}$ is said to posses a spectral gap (at $\JSR(\mathcal{A})$) 
if there exists $\gamma<1$ and
for every product $\tilde{\Pi}=\tilde{A}_{j_n}\cdots \tilde{A}_{j_1}$, 
$\tilde{A}_j\in\tilde{\mathcal{A}}$,
 which is not an s.m.p., 
it holds that
$\rho(\tilde{\Pi})<\gamma$.
\end{definition}


We are now in the position to describe the ipa,
which runs in two stages:
Firstly, it guesses spectral maximizing products $\Pi_n$, $n=1,\ldots,N$;
Secondly, it tries to construct the unit ball $P$ of a vector norm,
for whose induced matrix norm all normalized matrices $\tilde{A}_j=A_j/\rho(\Pi_1)^{1/\operatorname{len}(\Pi_1)}$, $A_j\in\mathcal{A}$, 
have norm less than or equal to $1$.
If the second parts succeeds, 
then, by Inequality~\eqref{equ_bound_jsr}, we obtain the exact value of the joint spectral radius.

The construction of the set $P$ is done iteratively. 
Starting with (properly scaled leading eigenvectors) of the s.m.p.-candidates, 
in each step it is checked whether all images of all points 
not yet mapped into the interior (of the convex hull of all formerly computed points)
are mapped into the interior
(of the convex hull of all formerly computed points).
Depending on the structure of the input set, different convex hulls need to be used;
\mbox{\emph{Case $(P)$}:~}If all entries of the matrices $A_j$ are non-negative, 
then we can take non-negative leading eigenvectors of the s.m.p.-candidates
as starting vectors and use the cone hull $\cocone$.
\mbox{\emph{Case $(R)$}:}~If the matrices $A_j$ have positive and negative entries 
and all leading eigenvectors are real, then we use the symmetric convex hull $\cosym$.
\mbox{\emph{Case $(C)$}:}~ In all other cases we need to use the elliptic convex hull $\coell$.
If eventually all points are mapped into the interior, 
then an invariant polytope is found and the algorithm terminates.

A simplified pseudo code implementation is given in Algorithm~\ref{alg_ipa}.
For a more thorough discussion of the algorithm see~\cite{GP2013,GP2016,Mej2020};
For a discussion about the containment problem see~\cite{GP2013,MP2022}.
\begin{algorithm}[!ht]
\caption{Invariant polytope algorithm\XX{checked}}
\label{alg_ipa}
\DontPrintSemicolon
\KwData{irreducible, finite set of matrices $\mathcal{A}=\set{A_j\in\RR^{s\times s}:j=1,\ldots,J}$}
\KwResultT{$\lambda=\JSR(\mathcal{A})$, invariant polytope $\coast V$}
Search for s.m.p.s $\Pi_1,\ldots, \Pi_N$, set $\lambda \coloneqq \rho(\Pi_1)^{1/\len\Pi_1}$\;
Scale matrices $\tilde{\mathcal{A}} \coloneqq \{\lambda^{-1} A_j:j=1,\ldots,J\}$\;
Select leading eigenvectors $V \coloneqq \set{v_0,\ldots,v_N}$\;
Set $R_{new} \coloneqq V$\;
\While{ $R_{new}\neq\emptyset$ }{
    Set $R \coloneqq R_{new}$\;
    Set $R_{new} \coloneqq \emptyset$\;
    \For{ $r\in\tilde{\mathcal{A}}R$ }{
        \If{ $r\notin\coast V^\circ$ }{
             Set $V \coloneqq V\cup r$\;
             Set $R_{new} \coloneqq R_{new} \cup r $\;
        }
    }
}
{\bfseries Return} $\lambda$, $\coast V$\;

\vspace{-5mm}
\end{algorithm}

\noindent
A crucial point in Algorithm~\ref{alg_ipa} is the line
\mbox{%
{\bfseries if} 
$r \notin \coast V^\circ$
{\bfseries then}
}%
:
If we cannot proof that a vector $r$
is not contained in the \emph{interior} of $\coast V$, then we have to add it to the set $V$.
Otherwise we would not get rigorous results using this algorithm.
Furthermore, with this procedure
sufficient and necessary conditions for the termination of the ipa are known.
\begin{theorem}[\cite{GP2016}]
\label{thm_termination}
Let $\mathcal{A}=\set{A_j\in\RR^{s\times s}:j=1,\ldots,J}$ be a finite set of matrices.
The ipa terminates if and only if 
the set $\mathcal{A}$ 
\begin{itemize}
\item has a spectral gap,
\item has only finitely many s.m.p.s $\Pi_n$, $n=1,\ldots,N$ (up to powers and cyclic permutations), and
\item each s.m.p.\ has only one simple leading eigenvector $v_n$ (up to complex conjugates).
\end{itemize}
\end{theorem}

\subsection{Mixed numeric/symbolic computations}
\label{sec_symbolic}
The conditions on the matrix set $\mathcal{A}$ in Theorem~\ref{thm_termination} sound rather restricting,
but it turns out that most matrix families from applications fulfil them.
Notable exceptions are when the scaled set $\tilde{\mathcal{A}}$ has a matrix product which is the identity matrix,
or
when vertices are mapped onto the boundary of the current polytope.
In both cases the ipa cannot terminate, since the algorithm always checks 
whether images of vertices are mapped into the interior of the current polytope.

To overcome these problems one can revert to a symbolic computation of the norm.
Unfortunately, a purely symbolic computation is computationally not feasible, 
because too expensive.
Thus, we resort to a mixed numerical and symbolic algorithm 
to replace the aforementioned line in the algorithm with
\mbox{%
{\bfseries if}
$s \notin \coast V$
{\bfseries then}%
} whenever possible.
We distinguish between two cases.

Case 1: Whenever a new vertex point is near to an existing vertex point,
we compare their exact coordinates symbolically. 
This can be done efficiently and just needs some matrix-vector multiplications.

Case 2: Slightly more complicated but still feasible;
Whenever the norm of a new vertex is near~$1$, 
we compute an exact upper bound of its norm symbolically.
This is efficiently possible whenever the leading eigenvectors are all real,
i.e.\ in cases $(R)$ and $(P)$.
Indeed, the problem of determining whether a point is inside or outside of a polytope 
can be stated as an LP problem~\cite{GP2013},
which does not only answer the containment problem, 
but also reports the vertices of a face of the polytope through 
which a ray through the point in question passes.
This face can be used to symbolically compute an \emph{upper} bound of the norm.
For the case when one leading eigenvector is complex (case $(C)$)
we yet do not have devised an efficient algorithm for the second problem.

Example~\ref{ex_prod_unitmatrix} shows how mixed numeric/symbolic computation 
can be used to solve examples where vertices of the polytope are mapped onto other vertices.
\begin{example}
\label{ex_prod_unitmatrix}
\XX{checked}
Let 
\begin{equation*}
A_1=\frac{1}{\sqrt{2}}\begin{bmatrix}0&0&-1\\0&0& 0\\ 0&1& 0\end{bmatrix},\quad
A_2=\frac{1}{\sqrt{2}}\begin{bmatrix}1&0&-1\\0&0&-1\\-1&0&-1\end{bmatrix}.
\end{equation*}
The only s.m.p.\ of this set is given by $A_2$, see below for the proof.
The ipa cannot compute the joint spectral radius of this set exactly due to two reasons:
$(1)$ The matrix $A_2$ has multiple leading eigenvalues $\pm1$, and furthermore
$(2)$, $A_2=A_2^3$ and thus vertices of any polytope are mapped onto itself after three iterations. 

With mixed symbolic and numeric computation we obtain that,
with leading eigenvector
$v_0=\begin{bmatrix}\sqrt{2}+2&\sqrt{2}&-2\end{bmatrix}^T$,
the polytope
$P=\cosym\set{
v_0, A_1 v_0, A_1 A_1 v_0, A_2 A_1 v_0, A_2 A_2 A_1 v_0
}$
is invariant under both matrices $A_1$, $A_2$. 
See Figure~\ref{fig_prod_unitmatrix} for the tree generated by the ipa.
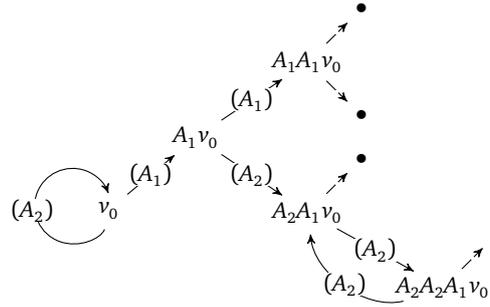
\begin{figure}[htb]
\centering
\begin{tikzpicture}[>=stealth',node distance=10mm, inner sep=0pt,
        ]

    \tikzset{mynode/.style={rectangle,inner sep=0pt}}
    \tikzset{mylastnode/.style={circle,inner sep=0pt}}    
    \node (T0)  [mynode] {$v_0$};
    \node (T1)  [mynode, above right=of T0] {$A_1 v_0$};
    \node (T11)  [mynode, above right=of T1] {$A_1A_1 v_0$};
    \node (T111) [mynode, above right of= T11] {$\bullet$};
    \node (T112) [mynode, below right of= T11] {$\bullet$};
    \node (T12) [mynode, below right=of T1] {$A_2A_1 v_0$};
    \node (T121) [mynode, above right of= T12] {$\bullet$};
    \node (T122) [mynode, below right=of T12] {$A_2 A_2 A_1 v_0$};
    \node (T1221) [mynode, above right of= T122] {$\bullet$};

    \tikzset{mypath/.style={->,shorten >=1ex, shorten <=1ex}}
    \tikzset{every node/.style={fill=white}}
    \draw[->,shorten <=1ex, shorten >=0ex] (T0.270) arc (340:20:5mm);
    \node[mynode] at (-1,0) {$(A_2)$};
    \path (T0) edge [mypath] node {$(A_1)$} (T1);
    \path (T1) edge [mypath] node {$(A_1)$} (T11);
    \path (T11) edge [mypath] node {} (T111);
    \path (T11) edge [mypath] node {} (T112);
    \path (T1) edge [mypath] node {$(A_2)$}  (T12);
    \path (T12) edge [mypath] node {$(A_2)$}  (T122);
    \path (T12) edge [mypath] node {}  (T121);
    \path (T122) edge [mypath,bend left=50pt] node {$(A_2)$}  (T12);
    \path (T122) edge [mypath] node {}  (T1221);
\end{tikzpicture}
\caption{Tree generated by the ipa with mixed numeric/symbolic computations for Example~\ref{ex_prod_unitmatrix}.
The starting vector $v_0$ is the leading eigenvector of $A_2$.
Arrows depict how vertices are mapped under the given matrix product.
Vertices plotted as $\bullet$ (instead written as text), are mapped to the interior of the polytope
$P=\cosym\set{
v_0, A_1 v_0, A_1 A_1 v_0, A_2 A_1 v_0, A_2 A_2 A_1 v_0
}$.
}
\label{fig_prod_unitmatrix}
\end{figure}
\end{example}
\begin{proof}
\label{proof_prod_unitmatrix}
We prove that $A_2$ is the only s.m.p.~of the set $\set{A_1,A_2}$.
First note that $\rho(A_2)=1$.
The claim follows by Gripenberg's algorithm~\cite{Grip1996}:
Since the norm $\norm{A_1}_2=\sqrt{2}/2<\rho(A_2)$, each product which is a candidate for an s.m.p.
has to start with either $\cdots A_1 A_2$ or $\cdots A_2 A_2$, where
\begin{align*}
A_1 A_2^1=\frac{1}{2}\begin{bmatrix}1&0&1\\0&0& 0\\ 0&0&-1\end{bmatrix},\quad \text{and} \quad
A_2 A_2^1=\frac{1}{2}\begin{bmatrix}2&0&0\\1&0&1\\0&0&2\end{bmatrix}.
\end{align*}
Again, $\norm{A_1 A_2^1}_2 = \sqrt{\frac{\sqrt{5} + 3}{8}}\simeq 0.80902 < 1$, and thus 
each product which is a candidate for an s.m.p.
has to start with either $\cdots A_1 A_2^2$ or $\cdots A_2 A_2^2$, where
\begin{align*}
A_1 A_2^2=\frac{1}{2\sqrt{2}}\begin{bmatrix}0&0&-2\\0&0& 0\\ 1&0&1\end{bmatrix},\quad \text{and} \quad
A_2 A_2^2=\frac{1}{\sqrt{2}}\begin{bmatrix}1&0&-1\\0&0&-1\\-1&0&-1\end{bmatrix}.
\end{align*}
And again, the norm $\norm{A_1 A_2^2}_2 =  \sqrt{\frac{\sqrt{5} + 3}{8}}\simeq 0.80902 < 1$, and thus
each product which is a candidate for an s.m.p.
has to start with either $\cdots A_1 A_2^3$ or $\cdots A_2 A_2^3$.
Since $A_1 A_2^3 = A_1 A_2^1$ and $A_2 A_2^3 = A_2 A_2^1$ we conclude that all s.m.p. candidates are of the form
$A_1 A_2^n$ and $A_2 A_2^n$.
The former are no s.m.p.s, since their spectral radii is $1/2$,
thus $A_2$ is the only s.m.p..
\end{proof}

\subsection{Limit matrices}
\label{sec_limitmatrix}

In some cases it speeds up the computation when one adds matrices 
to the input set $\mathcal{A}$ in question.
In particular, given a set of matrices $\mathcal{A}$, 
its joint spectral radius $\JSR(\mathcal{A})$ does not change
when elements of the closure $\closure \mathcal{A}$ 
or its convex hull $\co\mathcal{A}$
(to be understood in the Hausdorff distance using a matrix norm)
are added to $\mathcal{A}$~\cite[Proposition 1.8]{Jung2009},
\begin{equation}
\label{equ_closure}
\JSR(\mathcal{A}) =
\JSR(\closure\mathcal{A}) =
\JSR(\co\mathcal{A}).
\end{equation}
Lemma~\ref{thm_ipalimit} shows that we may also add limit matrices to the set in question.

\begin{lemma}
\label{thm_ipalimit}
Given matrices $\tilde{\mathcal{A}}\subseteq\RR^{s\times s}$ with 
$\JSR(\tilde{\mathcal{A}})=1$;
If the ipa terminates for the set $\tilde{\mathcal{A}}$, 
then the ipa terminates for the set
$\mathcal{L}\cup\tilde{\mathcal{A}}$,
where $\mathcal{L}$ is the set of all matrices of the limit set of 
all possible products of matrices of $\tilde{\mathcal{A}}$, 
i.e.\ of the set
$\set{ \prod_{k=1}^{n} \tilde{A}_{j_k} : 
\allowbreak n \in\NN,\ \allowbreak \tilde{A}_j\in\tilde{\mathcal{A}} }$.
\end{lemma}
\begin{proof}
\XX{checked}
If the ipa terminates, then there exists $K\in\NN$ such that 
$\tilde{A}_{j} v_n\in \coast V$ for all $\tilde{A}_j\in\tilde{\mathcal{A}}$ with
$
V=
\bigcup_{k=0}^{K} \tilde{\mathcal{A}}^k \{v_1,\ldots,v_N\}
$,
where $v_n$, $n=1,\ldots,N$, are the starting vectors for the ipa.
Thus, for each $v_n\in V$, there exists $M\in\NN$ such that
$\tilde{\Pi}_j^m v\in \coast  V$ for all $m \geq M$
and all s.m.p.s $\tilde{\Pi}_j$.
In particular, $\tilde{\Pi}_l v\in \coast V$ for all $\tilde{\Pi}_l \in \mathcal{L}$.
\end{proof}

Currently we use a heuristic to decide whether to add matrices and which of them.
The according rules have not stabilized yet, and thus, we do not report them.


\section{The finiteness conjecture}
\label{sec_fc}

Recalling the definition of a spectral maximizing product (s.m.p.) in Section~\ref{sec_ipa},
we make the following definition:
\begin{definition}
\label{def_fp}
A bounded set of matrices $\mathcal{A}\subseteq\RR^{s\times s}$ is said to posses the 
\emph{finiteness property} if there exists a finite product 
$\Pi=A_{j_n}\cdots A_{j_1}$, $A_{j_i}\in\mathcal{A}$ such that $\rho(\Pi)^{1/n}=\JSR(\mathcal{A})$.
\end{definition}

As already mentioned in the beginning, it has been shown that 
there exist sets of matrices such that the normalized spectral radius of every 
finite product is strictly less than the $\JSR$.
In other words, not all sets of matrices posses an s.m.p..
It is an open question whether pairs of \emph{binary matrices} with entries $\{0,1\}$ 
or \emph{sign matrices} with entries $\{-1,0,1\}$ always posses an s.m.p.~\cite{JB2008,Jung2009}.
Using the ipa we can check special cases of this question.

\begin{theorem}
\label{thm_fc}
The finiteness conjecture holds for all pairs of
\begin{itemize}
\item[$(a)$] binary matrices of dimension 2 $($i.e.\ with entries $\set{0, 1}$$)$,
\item[$(b)$] sign matrices of dimension 2 $($i.e.\ with entries $\set{-1, 0, 1}$$)$, and
\item[$(c)$] binary matrices of dimension 3 $($i.e.\ with entries $\set{0, 1}$$)$.
\end{itemize}
\end{theorem}
\begin{remark}
Point~\ref{thm_fc}~$(a)$ is already proven in~\cite[Chapter 4]{Jung2009};
Point~\ref{thm_fc}~$(b)$ is already proven in~\cite{CCGZ2010}.
\end{remark}

\begin{proof}
With our proposed algorithm a proof of $(a)$ and $(b)$ takes some minutes.
The proof of~\ref{thm_fc}~$(c)$ takes two days (CPU: AMD Ryzen 3600, 6 cores, 64 GB RAM) .
The used scripts to proof the results can be found at 
\href{https://gitlab.com/tommsch/dolomites}{\emph{gitlab.com/tommsch/dolomites}}
(and/or 
\href{http://tommsch.com/science.php}{\emph{tommsch.com/science.php}}%
)
and are named
\verb|fc_2.m|,
\verb|fc_2s.m|,
\verb|fc_3.m|.
The results in condensed form are tabulated in the Appendix;
in more detail they can be found online.
\end{proof}

\begin{remark}
If the algorithm would be implemented in a performant language (like C), 
this approach of checking the finiteness conjecture could also be used for 
pairs of sign matrices of dimension 3,
of which there are approximately 20 million cases to be checked. 
For larger matrices, this approach is not feasible any more.
\end{remark}

\subsection{Diminishing the number of cases}
\label{sec_diminish}
To proof~\ref{thm_fc}~$(c)$ we have to consider $2^{18}=262144$ cases
(To proof~$(b)$ we have to consider $3^8=6561$ cases, for $(a)$ $2^8=256$ cases).
This number can be reduced significantly: 
For some sets of matrices a concrete s.m.p.\ is known,
other sets share certain symmetries,
so that in total we check 15908 cases
(For $(b)$ we check $166$ cases, for $(a)$ we check $6$ cases).
We could exploit even more symmetries, but since those are computational hard to check,
the total time needed to proof the statement most likely would increase.

\begin{lemma}
\label{thm_diminish_A}
\XX{checked}
Given $A_1,A_2\in\RR^{s\times s}$; The following pairs have the same joint spectral radius
and the finiteness property holds for all or none of them:
\setlength{\columnsep}{-10cm}
\begin{multicols}{2}
\begin{itemize}
\item $\set{A_1,A_2}$,
\item $\set{A_2,A_1}$,
\item $\set{A_1^T,A_2^T}$,
\item $\set{\pm A_1,\pm A_2}$,
\item $\set{P^TA_1P,P^TA_2P}$ where $P$ is a permutation matrix, and
\item $\set{S^{-1}A_1S,S^{-1}A_2S}$ where $S$ is an invertible matrix.
\end{itemize}
\end{multicols}
\end{lemma}
\begin{proof}
For the proof we use
the definition of the joint spectral radius~\ref{def_jsr}
with the $2$-norm. The statements then follow from the facts that
$\norm{A}_2 = \norm{A^T}_2 = \norm{-A}_2$
and $P P^T = S S^{-1} = I$.
\end{proof}

\begin{lemma}
\label{thm_diminish_E}
\XX{checked}
\label{thm_jsr_monotone}
Given $A_1,A_2,A_0\in\NN_0^{s\times s}$;
If $A_2\leq A_1$, then $\JSR(\set{A_2,A_0})\leq \JSR(\set{A_1,A_0})$.
\end{lemma}
\begin{proof}
For the proof we use
the definition of the joint spectral radius~\ref{def_jsr}
with the Frobenius norm $\norm{\vardot}_F$,
and let $X = A_{j_n}\cdots A_{j_1}$, $j_i\in\set{2,0}$, be a given product.
We first construct a new product $\tilde{X}  = A_{\tilde{j}_n}\cdots A_{\tilde{j}_1}$, $\tilde{j}_i\in\set{1,0}$, from $X$, 
by replacing all occurrences of $A_2$ by $A_1$.
It follows that
$
\norm{X}_F^{1/n} \leq \norm{\tilde{X}}_F^{1/n}$,
and thus $\JSR(\set{A_2,A_0}) \leq \JSR(\set{A_1,A_0})$.
\end{proof}

\begin{lemma}
\XX{checked}
Given $A_1,A_2\in\NN_0^{s\times s}$;
The finiteness property holds whenever $\JSR(\{A_1,A_2\})\leq 1$.
\end{lemma}
\begin{proof}
Since the norm of a non-zero integer matrix is always greater equal than one, 
it is not possible that the joint spectral radius of a set of integer matrices is strictly between $0$ and $1$.
If $\JSR(\{A_1,A_2\})=0$, then clearly both $A_1$ and $A_2$ are s.m.p.s.
The second case is non trivial and its proof is given in~\cite[Chapter 3.4]{Jung2009}.
\end{proof}

\begin{corollary}
\XX{checked}
Given $A_1,A_2\in\NN_0^{s\times s}$;
The finiteness property holds whenever
\setlength{\columnsep}{-10cm}
\begin{multicols}{2}
\begin{itemize}
\item[$(a)$] $A_2\leq A_1$
\item[$(b)$] $A_1A_2\leq A_1^2$ 
\item[$(c)$] $A_2\leq I$
\item[$(d)$] $A_2A_1\leq A_1A_2$ 
\end{itemize}
\end{multicols}
\end{corollary}
\begin{proof}
Again, we use
the definition of the joint spectral radius~\ref{def_jsr}
with the Frobenius norm $\norm{\vardot}_F$,
and let $A_{j_n}\cdots A_{j_1}$, $j\in\set{1,2}$, be a given product.

$(a)$ and $(b)$ It follows that
$
\norm{A_{j_n}\cdots A_{j_1}}_F^{1/n} \leq 
\norm{A_1}_F^{1/n},
$
and thus $\JSR(\set{A_1,A_2}) = \rho(A_1)$

$(c)$ It follows that
$
\norm{A_{j_n}\cdots A_{j_1}}_F^{1/n} \leq 
\norm{A_1}_F^{1/\tilde{n}}
$
with $\tilde{n}\leq n$,
and thus, $\JSR(\set{A_1,A_2})\leq \rho(A_1)$ which implies
$\JSR(\set{A_1,A_2}) = \rho(A_1)$.

$(d)$ It follows that
$
\norm{A_{j_n}\cdots A_{j_1}}_F^{1/n} \leq 
\norm{A_1^{n_1} A_2^{n_2}}_F^{1/n}
$
for some $n_1 + n_2 = n$,
and thus,
$\JSR(\set{A_1,A_2}) = \max\set{\rho(A_1),\rho(A_2)}$.
\end{proof}

\begin{lemma}
\label{thm_norm_fp}
\XX{checked}
If there exists a norm $\norm{\vardot}$ such that 
$\max \set{\rho(A_1),\rho(A_2)} = \max \set{\norm{A_1},\norm{A_2}}$, then the finiteness property holds.
In particular, the finiteness property holds for sets of
normal matrices, and thus, symmetric matrices.
A matrix $A$ is normal, iff $A^TA = AA^T$.
\end{lemma}
\begin{proof}
The first part follows from Inequality~\eqref{equ_bound_jsr},
which reads for products of length~1 as
$
\max_{A_{j}\in\mathcal{A}}\rho\left( A_j\right)
\leq
\JSR(\mathcal{A})
\leq
\max_{A_j\in\mathcal{A}}
\norm{A_{j}}
$.
By the assumptions we have equality here, and thus
$\JSR(\set{A_1,A_2}) = \max_{A_{j}\in\mathcal{A}}\rho\left( A_j\right)$.

The second parts about normal matrices follows now by using the $2$-norm,
which equals the matrix' the largest singular value.
For normal matrices the largest singular value equals the largest eigenvalue in magnitude,
and thus
$
\max_{A_{j}\in\mathcal{A}}\rho\left( A_j\right) = 
\max_{A_j\in\mathcal{A}} \norm{ A_j }_2$.
\end{proof}

\begin{definition}
\label{def_irr}
Given a finite set of matrices $\mathcal{A}\subseteq\RR^{s\times s}$;
If there exists $V\in\RR^{s\times s}$ such that
$V A_j V^{-1} = \begin{bmatrix} B_j & C_j\\0 & D_j\end{bmatrix}$ for all $A_j\in\mathcal{A}$,
then $\mathcal{A}$ is \emph{reducible}.
\end{definition}

\begin{theorem}
\label{thm_diminish_T}
In the notation from Definition~\ref{def_irr};
If $\mathcal{A}$ is reducible, then 
$\JSR(\mathcal{A}) = \max \{\JSR(\mathcal{B}),\allowbreak \JSR(\mathcal{D})\}$,
$\mathcal{B} = \set{ B_j : j = 1,\ldots,J}$,
$\mathcal{D} = \set{ D_j : j = 1,\ldots, J}$.
\end{theorem}
The first rigorous proof known to the author can be found in~\cite[Proposition 1.5]{Jung2009}.
Although the proof is straight forward, it is also rather technical and we abstain from giving it here.

\begin{corollary}
\XX{checked}
Given $A_1,A_2\in Z^{s\times s}$, $Z\subseteq \ZZ$
and using the notation from Definition~\ref{def_irr};
The finiteness conjecture holds whenever there exists
$S\in\CC^{s\times s}$ such that the matrices
$S^{-1}A_jS$, have joint block diagonal form with blocks 
$B_j\in Z^{s_B \times s_B}$, 
$D_j\in Z^{s_D \times s_D}$, $j=1,2$, $s_B<l$,
and the finites property holds for all pairs of matrices in $Z^{s_B\times s_B}$.
\end{corollary}
\begin{proof}
This follows from Lemma~\ref{thm_diminish_A} and Theorem~\ref{thm_diminish_T}.
\end{proof}

\section{Implementation notes}
Our Matlab implementation of the algorithm can be found on Gitlab~\cite{ttoolboxes},
and is extensively documented.
The file \texttt{manual.pdf} gives an overview of the toolbox,
in depth documentation can be found directly in the source files,
and can be viewed by typing \texttt{help functionname} or \texttt{edit functionname},
e.g.\ \texttt{help tjsr} or \texttt{edit tjsr}, in Matlab.

The main function for the JSR computation is~\texttt{tjsr}, 
short for \emph{invarianT polyTope algoriThm}.
Depending on the input, our implementation chooses its parameters automatically 
and usually there is no need for the user to specify options by hand.
For example, of the 15910 cases checked for Theorem~\ref{thm_fc}~$(c)$,
manual intervention was only necessary for 2 cases.

\begin{figure}
\centering
\includegraphics[scale=0.8]{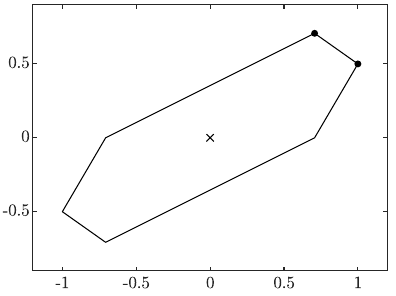}
\caption%
[Invariant polytope produced by the ipa for the matrices in Example~\ref{ex_tjsr}]%
{Invariant polytope for the matrices 
$A_1 = 
\begin{bsmallmatrix}
0 & 1\\0 & 1
\end{bsmallmatrix},\quad
A_2 = \begin{bsmallmatrix}
1&\phantom{-}0\\1&-1
\end{bsmallmatrix}.
$
of Example~\ref{ex_tjsr}.
The cross~$\times$ denotes the origin, the dots $\bullet$ the leading eigenvectors of the matrices
(and s.m.p.s)
$A_1$ and~$A_2$.}
\label{fig_tjsr}
\vspace{-5mm}
\end{figure}

Example~\ref{ex_tjsr} presents how to use the \texttt{tjsr} algorithm.
\begin{example}
\label{ex_tjsr}
Given the matrices 
$
A_1 = 
\begin{bsmallmatrix}
0 & 1\\0 & 1
\end{bsmallmatrix},\quad
A_2 = \begin{bsmallmatrix}
1&\phantom{-}0\\1&-1
\end{bsmallmatrix}.
$
To compute their joint spectral radius, 
the matrices must be passed as a cell array to the algorithm,
e.g.\ by typing:
{\small%
\begin{verbatim}
    tjsr( {[0 1;0 1],[1 0;1 -1]} )
\end{verbatim}%
}%
\noindent
The algorithm (version 1.2022.05.25) produces the following output:

\definecolor{dred}{rgb}{0.7, 0.0, 0.0}
\definecolor{dgreen}{rgb}{0.0, 0.5, 0.0}
\definecolor{dorange}{rgb}{0.8, 0.4, 0.3}
{\footnotesize
\begin{Verbatim}[commandchars=\\\{\}]
    Input: 2 matrices of dimension 2
    
    {\color{dred}A lot of candidates found. Nearly all orderings are smp's.}
    {\color{dred}Set <'epseigenplane',inf, 'epsspectralradius',inf, 'maxsmpdepth',5, 'balancing',-1, 'ncdelta',1>.}
    JSR (of block) 1:     1.0000    1.6180
    {\color{dred}Duplicate leading eigenvectors occured. Enable symbolic computation. Set <'epssym',5e-12>.}
    JSR (of block) 1:     1.0000    1.6180
    {\color{dorange}Case (R).}
    Selected candidates: | 1 | 2 | 
    Number of vertices: 2 ( candidates )
    Balance 2 Trees. Balancing vector found: [ 1, 341/305]
    JSR = [               1,   1.61803398875 ], norm=           Inf, #test: 1/1, #V:2/2 | O
    JSR = [               1,   1.61803398875 ], norm= 2.41421508763, _
    Number of vertices of polytope: 3
    Products which give lower bounds of JSR: | 1 |
    {\color{dgreen}Algorithm terminated correctly. Exact value found.}
    JSR =               1
\end{Verbatim}
}

One can see that the algorithm restarts two times.
The first time because two many s.m.p.\ candidates are found
and appropriate options are set:
\texttt{\small{\color{dred}<%
\textquotesingle epseigenplane\textquotesingle,inf, 
\textquotesingle epsspectralradius\textquotesingle ,inf, 
\textquotesingle maxsmpdepth\textquotesingle ,5, 
\textquotesingle balancing\textquotesingle ,-1, 
\textquotesingle ncdelta\textquotesingle ,1>}}.
The second time because the s.m.p.\ candidates do not have a unique leading eigenvector
and mixed symbolic and numeric computation is enabled
{\texttt{\small\color{dred}<%
\textquotesingle epssym\textquotesingle,5e-12>}}.

The third time the algorithm terminates after two iterations.
It reports two s.m.p.s $A_1$ and $A_2$, and a joint spectral radius of $1$.
The constructed polytope with 3 vertices is given in Figure~\ref{fig_tjsr}.
The figure is produced by calling \texttt{tjsr} with the option \texttt{\textquotesingle plot\textquotesingle}:
\texttt{tjsr( {[0 1;0 1],[1 0;1 -1]}, \textquotesingle plot\textquotesingle,\textquotesingle polytope\textquotesingle{} )}. A complete list of all options
can be found in the file \texttt{tjsr\_option},
and viewed by typing \texttt{tjsr help} or \texttt{edit tjsr\_option}.
\end{example}

\begin{remark}
Since our implementation of the algorithm is in Matlab, which has very restricted capabilities
for symbolic computations, mixed symbolic computation only works
when the leading eigenvectors are expressible
in a ``simple'' closed form, 
e.g.\ for integer matrices of dimension less than or equal to~3.
\end{remark}

\appendix

\section{List of cases}

\subsection{\texorpdfstring
    {$2\times2$ binary matrices}%
    {2 times 2 binary matrices}%
}
\label{sec_f2s}
The following list reports 
$\mathcal{A}$: the set of matrices and
$s.m.p.$: the shortest s.m.p. found for the set $\mathcal{A}$.
All unreported cases can be reduced to a simpler one,
or an s.m.p. is known due to the structure of the set $\mathcal{A}$,
by the Lemmata presented in Section~\ref{sec_diminish}.

{%
\centering
\small
     \begin{longtable}{lc | lc | lc }
     \hfil$\mathcal{A}$\hfil & \emph{s.m.p.}\Bstrut & 
        \hfil$\mathcal{A}$\hfil  & \emph{s.m.p.}\Bstrut &
            \hfil$\mathcal{A}$\hfil & \emph{s.m.p.}\Bstrut \\ \hline
     $\{\begin{bsmallmatrix}0&1\\0&0\end{bsmallmatrix},\ \begin{bsmallmatrix}1&0\\1&1\end{bsmallmatrix}\}$ & $A_{1}A_{2}^{4}$\Tstrut &
         $\{\begin{bsmallmatrix}0&1\\1&0\end{bsmallmatrix},\ \begin{bsmallmatrix}1&1\\0&1\end{bsmallmatrix}\}$ & $A_{1}A_{2}^{3}$ &
                $\{\begin{bsmallmatrix}0&1\\0&1\end{bsmallmatrix},\ \begin{bsmallmatrix}0&0\\1&1\end{bsmallmatrix}\}$ & $A_{1}A_{2}$ \\
  $\{\begin{bsmallmatrix}1&1\\0&1\end{bsmallmatrix},\ \begin{bsmallmatrix}1&0\\1&0\end{bsmallmatrix}\}$  & $A_{1}^{2}A_{2}$ &
     $\{\begin{bsmallmatrix}0&1\\0&1\end{bsmallmatrix},\ \begin{bsmallmatrix}1&0\\1&1\end{bsmallmatrix}\}$  & $A_{1}A_{2}^{2}$ &
        $\{\begin{bsmallmatrix}1&1\\0&1\end{bsmallmatrix},\ \begin{bsmallmatrix}1&0\\1&1\end{bsmallmatrix}\}$& $A_{1}A_{2}$ \\
     \end{longtable}

\par
}

\subsection{\texorpdfstring
    {$2\times2$ sign matrices}%
    {2 times 2 sign matrices}%
}
To save space, we use the following abbreviations: $+$ for $+1$, $-$ for $-1$, $\circ$ for $0$.
In addition to $\mathcal{A}$ and $s.m.p.$, as reported in Section~\ref{sec_f2s},
we mark all cases where at least one leading eigenvalues is complex (case $(C)$).
If the case is not explicitly mentioned, then it is case $(R)$, meaning that all leading eigenvalues are real.

\bigskip%
{%
\centering%
\tiny%

\newcommand\pl{\makebox[0.8em]{$+$}}
\newcommand\mn{\makebox[0.8em]{$-$}}
\newcommand\zr{\makebox[0.8em]{$\circ$}}
\begin{xtabular}{lcc|lcc|lcc|lcc|lcc}
\hfil $\mathcal{A}$\hfil & \!\!\!\emph{s.m.p.}\!\!\! & \!\!\!\!\emph{case}\!\!\!\! &\hfil $\mathcal{A}$\hfil & \!\!\!\emph{s.m.p.}\!\!\! & \!\!\!\!\emph{case}\!\!\!\! &\hfil $\mathcal{A}$\hfil & \!\!\!\emph{s.m.p.}\!\!\! & \!\!\!\!\emph{case}\!\!\!\! &\hfil $\mathcal{A}$\hfil & \!\!\!\emph{s.m.p.}\!\!\! & \!\!\!\!\emph{case}\!\!\!\! &\hfil $\mathcal{A}$\hfil & \!\!\!\emph{s.m.p.}\!\!\! & \!\!\!\!\emph{case}\!\!\!\! \Bstrut \\ \hline
$\!\!\!\{\!\begin{bsmallmatrix}\!\zr&\zr\!\\\!\zr&\pl\!\end{bsmallmatrix},\!\begin{bsmallmatrix}\!\zr&\mn\!\\\!\pl&\pl\!\end{bsmallmatrix}\!\}\!\!\!\!\!\!\!\!$ & $A_{1}\!\!\!$ & \!\!(C)\!\!&$\!\!\!\{\!\begin{bsmallmatrix}\!\zr&\zr\!\\\!\zr&\pl\!\end{bsmallmatrix},\!\begin{bsmallmatrix}\!\zr&\mn\!\\\!\pl&\mn\!\end{bsmallmatrix}\!\}\!\!\!\!\!\!\!\!$ & $A_{1}\!\!\!$ & \!\!(C)\!\!&$\!\!\!\{\!\begin{bsmallmatrix}\!\zr&\zr\!\\\!\zr&\pl\!\end{bsmallmatrix},\!\begin{bsmallmatrix}\!\pl&\mn\!\\\!\pl&\zr\!\end{bsmallmatrix}\!\}\!\!\!\!\!\!\!\!$ & $A_{1}\!\!\!$ & \!\!(C)\!\!&$\!\!\!\{\!\begin{bsmallmatrix}\!\zr&\zr\!\\\!\zr&\pl\!\end{bsmallmatrix},\!\begin{bsmallmatrix}\!\pl&\mn\!\\\!\pl&\mn\!\end{bsmallmatrix}\!\}\!\!\!\!\!\!\!\!$ & $A_{1}\!\!\!$ & &$\!\!\!\{\!\begin{bsmallmatrix}\!\zr&\pl\!\\\!\zr&\zr\!\end{bsmallmatrix},\!\begin{bsmallmatrix}\!\zr&\mn\!\\\!\pl&\pl\!\end{bsmallmatrix}\!\}\!\!\!\!\!\!\!\!$ & $A_{2}\!\!\!$ & \!\!(C)\!\!\\
$\!\!\!\{\!\begin{bsmallmatrix}\!\zr&\pl\!\\\!\zr&\zr\!\end{bsmallmatrix},\!\begin{bsmallmatrix}\!\zr&\mn\!\\\!\pl&\mn\!\end{bsmallmatrix}\!\}\!\!\!\!\!\!\!\!$ & $A_{2}\!\!\!$ & \!\!(C)\!\!&$\!\!\!\{\!\begin{bsmallmatrix}\!\zr&\pl\!\\\!\zr&\zr\!\end{bsmallmatrix},\!\begin{bsmallmatrix}\!\pl&\zr\!\\\!\pl&\mn\!\end{bsmallmatrix}\!\}\!\!\!\!\!\!\!\!$ & $A_{2}\!\!\!$ & &$\!\!\!\{\!\begin{bsmallmatrix}\!\zr&\pl\!\\\!\zr&\zr\!\end{bsmallmatrix},\!\begin{bsmallmatrix}\!\pl&\mn\!\\\!\pl&\zr\!\end{bsmallmatrix}\!\}\!\!\!\!\!\!\!\!$ & $A_{2}\!\!\!$ & \!\!(C)\!\!&$\!\!\!\{\!\begin{bsmallmatrix}\!\zr&\pl\!\\\!\zr&\zr\!\end{bsmallmatrix},\!\begin{bsmallmatrix}\!\pl&\mn\!\\\!\pl&\mn\!\end{bsmallmatrix}\!\}\!\!\!\!\!\!\!\!$ & $A_{1}A_{2}\!\!\!$ & &$\!\!\!\{\!\begin{bsmallmatrix}\!\zr&\pl\!\\\!\zr&\zr\!\end{bsmallmatrix},\!\begin{bsmallmatrix}\!\pl&\zr\!\\\!\mn&\pl\!\end{bsmallmatrix}\!\}\!\!\!\!\!\!\!\!$ & $A_{1}A_{2}^{4}\!\!\!$ & \\
$\!\!\!\{\!\begin{bsmallmatrix}\!\zr&\pl\!\\\!\zr&\zr\!\end{bsmallmatrix},\!\begin{bsmallmatrix}\!\pl&\zr\!\\\!\mn&\mn\!\end{bsmallmatrix}\!\}\!\!\!\!\!\!\!\!$ & $A_{2}\!\!\!$ & &$\!\!\!\{\!\begin{bsmallmatrix}\!\zr&\pl\!\\\!\zr&\zr\!\end{bsmallmatrix},\!\begin{bsmallmatrix}\!\pl&\pl\!\\\!\mn&\zr\!\end{bsmallmatrix}\!\}\!\!\!\!\!\!\!\!$ & $A_{2}\!\!\!$ & \!\!(C)\!\!&$\!\!\!\{\!\begin{bsmallmatrix}\!\zr&\pl\!\\\!\zr&\zr\!\end{bsmallmatrix},\!\begin{bsmallmatrix}\!\pl&\pl\!\\\!\mn&\mn\!\end{bsmallmatrix}\!\}\!\!\!\!\!\!\!\!$ & $A_{1}A_{2}\!\!\!$ & &$\!\!\!\{\!\begin{bsmallmatrix}\!\zr&\pl\!\\\!\zr&\pl\!\end{bsmallmatrix},\!\begin{bsmallmatrix}\!\zr&\mn\!\\\!\pl&\pl\!\end{bsmallmatrix}\!\}\!\!\!\!\!\!\!\!$ & $A_{1}A_{2}\!\!\!$ & &$\!\!\!\{\!\begin{bsmallmatrix}\!\zr&\pl\!\\\!\zr&\pl\!\end{bsmallmatrix},\!\begin{bsmallmatrix}\!\zr&\mn\!\\\!\pl&\mn\!\end{bsmallmatrix}\!\}\!\!\!\!\!\!\!\!$ & $A_{1}\!\!\!$ & \!\!(C)\!\!\\
$\!\!\!\{\!\begin{bsmallmatrix}\!\zr&\pl\!\\\!\zr&\pl\!\end{bsmallmatrix},\!\begin{bsmallmatrix}\!\pl&\zr\!\\\!\pl&\mn\!\end{bsmallmatrix}\!\}\!\!\!\!\!\!\!\!$ & $A_{1}\!\!\!$ & &$\!\!\!\{\!\begin{bsmallmatrix}\!\zr&\pl\!\\\!\zr&\pl\!\end{bsmallmatrix},\!\begin{bsmallmatrix}\!\pl&\mn\!\\\!\pl&\zr\!\end{bsmallmatrix}\!\}\!\!\!\!\!\!\!\!$ & $A_{1}\!\!\!$ & \!\!(C)\!\!&$\!\!\!\{\!\begin{bsmallmatrix}\!\zr&\pl\!\\\!\zr&\pl\!\end{bsmallmatrix},\!\begin{bsmallmatrix}\!\pl&\zr\!\\\!\mn&\pl\!\end{bsmallmatrix}\!\}\!\!\!\!\!\!\!\!$ & $A_{1}A_{2}^{5}\!\!\!$ & &$\!\!\!\{\!\begin{bsmallmatrix}\!\zr&\pl\!\\\!\zr&\pl\!\end{bsmallmatrix},\!\begin{bsmallmatrix}\!\pl&\zr\!\\\!\mn&\mn\!\end{bsmallmatrix}\!\}\!\!\!\!\!\!\!\!$ & $A_{1}A_{2}\!\!\!$ & &$\!\!\!\{\!\begin{bsmallmatrix}\!\zr&\pl\!\\\!\zr&\pl\!\end{bsmallmatrix},\!\begin{bsmallmatrix}\!\pl&\pl\!\\\!\mn&\zr\!\end{bsmallmatrix}\!\}\!\!\!\!\!\!\!\!$ & $A_{1}A_{2}^{2}\!\!\!$ & \\
$\!\!\!\{\!\begin{bsmallmatrix}\!\zr&\pl\!\\\!\zr&\pl\!\end{bsmallmatrix},\!\begin{bsmallmatrix}\!\pl&\pl\!\\\!\mn&\mn\!\end{bsmallmatrix}\!\}\!\!\!\!\!\!\!\!$ & $A_{1}A_{2}\!\!\!$ & &$\!\!\!\{\!\begin{bsmallmatrix}\!\zr&\pl\!\\\!\zr&\mn\!\end{bsmallmatrix},\!\begin{bsmallmatrix}\!\zr&\zr\!\\\!\pl&\mn\!\end{bsmallmatrix}\!\}\!\!\!\!\!\!\!\!$ & $A_{1}A_{2}\!\!\!$ & &$\!\!\!\{\!\begin{bsmallmatrix}\!\zr&\pl\!\\\!\zr&\mn\!\end{bsmallmatrix},\!\begin{bsmallmatrix}\!\zr&\mn\!\\\!\pl&\pl\!\end{bsmallmatrix}\!\}\!\!\!\!\!\!\!\!$ & $A_{1}\!\!\!$ & \!\!(C)\!\!&$\!\!\!\{\!\begin{bsmallmatrix}\!\zr&\pl\!\\\!\zr&\mn\!\end{bsmallmatrix},\!\begin{bsmallmatrix}\!\zr&\mn\!\\\!\pl&\mn\!\end{bsmallmatrix}\!\}\!\!\!\!\!\!\!\!$ & $A_{1}A_{2}\!\!\!$ & &$\!\!\!\{\!\begin{bsmallmatrix}\!\zr&\pl\!\\\!\zr&\mn\!\end{bsmallmatrix},\!\begin{bsmallmatrix}\!\pl&\zr\!\\\!\pl&\pl\!\end{bsmallmatrix}\!\}\!\!\!\!\!\!\!\!$ & $A_{1}A_{2}^{5}\!\!\!$ & \\
$\!\!\!\{\!\begin{bsmallmatrix}\!\zr&\pl\!\\\!\zr&\mn\!\end{bsmallmatrix},\!\begin{bsmallmatrix}\!\pl&\zr\!\\\!\pl&\mn\!\end{bsmallmatrix}\!\}\!\!\!\!\!\!\!\!$ & $A_{1}A_{2}\!\!\!$ & &$\!\!\!\{\!\begin{bsmallmatrix}\!\zr&\pl\!\\\!\zr&\mn\!\end{bsmallmatrix},\!\begin{bsmallmatrix}\!\pl&\mn\!\\\!\pl&\zr\!\end{bsmallmatrix}\!\}\!\!\!\!\!\!\!\!$ & $A_{1}A_{2}^{2}\!\!\!$ & &$\!\!\!\{\!\begin{bsmallmatrix}\!\zr&\pl\!\\\!\zr&\mn\!\end{bsmallmatrix},\!\begin{bsmallmatrix}\!\pl&\mn\!\\\!\pl&\mn\!\end{bsmallmatrix}\!\}\!\!\!\!\!\!\!\!$ & $A_{1}A_{2}\!\!\!$ & &$\!\!\!\{\!\begin{bsmallmatrix}\!\zr&\pl\!\\\!\zr&\mn\!\end{bsmallmatrix},\!\begin{bsmallmatrix}\!\pl&\zr\!\\\!\mn&\pl\!\end{bsmallmatrix}\!\}\!\!\!\!\!\!\!\!$ & $A_{1}A_{2}^{2}\!\!\!$ & &$\!\!\!\{\!\begin{bsmallmatrix}\!\zr&\pl\!\\\!\zr&\mn\!\end{bsmallmatrix},\!\begin{bsmallmatrix}\!\pl&\zr\!\\\!\mn&\mn\!\end{bsmallmatrix}\!\}\!\!\!\!\!\!\!\!$ & $A_{1}\!\!\!$ & \\
$\!\!\!\{\!\begin{bsmallmatrix}\!\zr&\pl\!\\\!\zr&\mn\!\end{bsmallmatrix},\!\begin{bsmallmatrix}\!\pl&\pl\!\\\!\mn&\zr\!\end{bsmallmatrix}\!\}\!\!\!\!\!\!\!\!$ & $A_{1}\!\!\!$ & \!\!(C)\!\!&$\!\!\!\{\!\begin{bsmallmatrix}\!\zr&\pl\!\\\!\pl&\zr\!\end{bsmallmatrix},\!\begin{bsmallmatrix}\!\zr&\mn\!\\\!\pl&\pl\!\end{bsmallmatrix}\!\}\!\!\!\!\!\!\!\!$ & $A_{1}\!\!\!$ & \!\!(C)\!\!&$\!\!\!\{\!\begin{bsmallmatrix}\!\zr&\pl\!\\\!\pl&\zr\!\end{bsmallmatrix},\!\begin{bsmallmatrix}\!\zr&\mn\!\\\!\pl&\mn\!\end{bsmallmatrix}\!\}\!\!\!\!\!\!\!\!$ & $A_{1}\!\!\!$ & \!\!(C)\!\!&$\!\!\!\{\!\begin{bsmallmatrix}\!\zr&\pl\!\\\!\pl&\zr\!\end{bsmallmatrix},\!\begin{bsmallmatrix}\!\pl&\pl\!\\\!\zr&\mn\!\end{bsmallmatrix}\!\}\!\!\!\!\!\!\!\!$ & $A_{1}\!\!\!$ & \!\!(C)\!\!&$\!\!\!\{\!\begin{bsmallmatrix}\!\zr&\pl\!\\\!\pl&\zr\!\end{bsmallmatrix},\!\begin{bsmallmatrix}\!\pl&\mn\!\\\!\zr&\pl\!\end{bsmallmatrix}\!\}\!\!\!\!\!\!\!\!$ & $A_{1}A_{2}^{3}\!\!\!$ & \\
$\!\!\!\{\!\begin{bsmallmatrix}\!\zr&\pl\!\\\!\pl&\zr\!\end{bsmallmatrix},\!\begin{bsmallmatrix}\!\pl&\mn\!\\\!\zr&\mn\!\end{bsmallmatrix}\!\}\!\!\!\!\!\!\!\!$ & $A_{1}\!\!\!$ & \!\!(C)\!\!&$\!\!\!\{\!\begin{bsmallmatrix}\!\zr&\pl\!\\\!\pl&\pl\!\end{bsmallmatrix},\!\begin{bsmallmatrix}\!\pl&\mn\!\\\!\pl&\mn\!\end{bsmallmatrix}\!\}\!\!\!\!\!\!\!\!$ & $A_{1}\!\!\!$ & &$\!\!\!\{\!\begin{bsmallmatrix}\!\zr&\pl\!\\\!\pl&\mn\!\end{bsmallmatrix},\!\begin{bsmallmatrix}\!\pl&\mn\!\\\!\pl&\mn\!\end{bsmallmatrix}\!\}\!\!\!\!\!\!\!\!$ & $A_{1}\!\!\!$ & &$\!\!\!\{\!\begin{bsmallmatrix}\!\zr&\mn\!\\\!\pl&\zr\!\end{bsmallmatrix},\!\begin{bsmallmatrix}\!\zr&\mn\!\\\!\pl&\pl\!\end{bsmallmatrix}\!\}\!\!\!\!\!\!\!\!$ & \multicolumn{4}{c}{$ A_{1}A_{2}A_{1}A_{2}^{2}$} & \\
$\!\!\!\{\!\begin{bsmallmatrix}\!\zr&\mn\!\\\!\pl&\zr\!\end{bsmallmatrix},\!\begin{bsmallmatrix}\!\zr&\mn\!\\\!\pl&\mn\!\end{bsmallmatrix}\!\}\!\!\!\!\!\!\!\!$ & \multicolumn{4}{c}{$ A_{1}A_{2}A_{1}A_{2}^{2}$} & &$\!\!\!\{\!\begin{bsmallmatrix}\!\zr&\mn\!\\\!\pl&\zr\!\end{bsmallmatrix},\!\begin{bsmallmatrix}\!\pl&\pl\!\\\!\zr&\pl\!\end{bsmallmatrix}\!\}\!\!\!\!\!\!\!\!$ & $A_{1}A_{2}^{4}\!\!\!$ & &$\!\!\!\{\!\begin{bsmallmatrix}\!\zr&\mn\!\\\!\pl&\zr\!\end{bsmallmatrix},\!\begin{bsmallmatrix}\!\pl&\pl\!\\\!\zr&\mn\!\end{bsmallmatrix}\!\}\!\!\!\!\!\!\!\!$ & $A_{1}A_{2}\!\!\!$ & &$\!\!\!\{\!\begin{bsmallmatrix}\!\zr&\mn\!\\\!\pl&\zr\!\end{bsmallmatrix},\!\begin{bsmallmatrix}\!\pl&\mn\!\\\!\zr&\pl\!\end{bsmallmatrix}\!\}\!\!\!\!\!\!\!\!$ & $A_{1}A_{2}^{4}\!\!\!$ & \\
$\!\!\!\{\!\begin{bsmallmatrix}\!\zr&\mn\!\\\!\pl&\zr\!\end{bsmallmatrix},\!\begin{bsmallmatrix}\!\pl&\mn\!\\\!\zr&\mn\!\end{bsmallmatrix}\!\}\!\!\!\!\!\!\!\!$ & $A_{1}A_{2}\!\!\!$ & &$\!\!\!\{\!\begin{bsmallmatrix}\!\zr&\mn\!\\\!\pl&\zr\!\end{bsmallmatrix},\!\begin{bsmallmatrix}\!\pl&\zr\!\\\!\pl&\pl\!\end{bsmallmatrix}\!\}\!\!\!\!\!\!\!\!$ & $A_{1}A_{2}^{4}\!\!\!$ & &$\!\!\!\{\!\begin{bsmallmatrix}\!\zr&\mn\!\\\!\pl&\zr\!\end{bsmallmatrix},\!\begin{bsmallmatrix}\!\pl&\zr\!\\\!\pl&\mn\!\end{bsmallmatrix}\!\}\!\!\!\!\!\!\!\!$ & $A_{1}A_{2}\!\!\!$ & &$\!\!\!\{\!\begin{bsmallmatrix}\!\zr&\mn\!\\\!\pl&\zr\!\end{bsmallmatrix},\!\begin{bsmallmatrix}\!\pl&\mn\!\\\!\pl&\zr\!\end{bsmallmatrix}\!\}\!\!\!\!\!\!\!\!$ & \multicolumn{4}{c}{$ A_{1}A_{2}A_{1}A_{2}^{2}$} & \\
$\!\!\!\{\!\begin{bsmallmatrix}\!\zr&\mn\!\\\!\pl&\zr\!\end{bsmallmatrix},\!\begin{bsmallmatrix}\!\pl&\mn\!\\\!\pl&\mn\!\end{bsmallmatrix}\!\}\!\!\!\!\!\!\!\!$ & $A_{1}A_{2}\!\!\!$ & &$\!\!\!\{\!\begin{bsmallmatrix}\!\zr&\mn\!\\\!\pl&\zr\!\end{bsmallmatrix},\!\begin{bsmallmatrix}\!\pl&\zr\!\\\!\mn&\pl\!\end{bsmallmatrix}\!\}\!\!\!\!\!\!\!\!$ & $A_{1}A_{2}^{4}\!\!\!$ & &$\!\!\!\{\!\begin{bsmallmatrix}\!\zr&\mn\!\\\!\pl&\zr\!\end{bsmallmatrix},\!\begin{bsmallmatrix}\!\pl&\zr\!\\\!\mn&\mn\!\end{bsmallmatrix}\!\}\!\!\!\!\!\!\!\!$ & $A_{1}A_{2}\!\!\!$ & &$\!\!\!\{\!\begin{bsmallmatrix}\!\zr&\mn\!\\\!\pl&\zr\!\end{bsmallmatrix},\!\begin{bsmallmatrix}\!\pl&\pl\!\\\!\mn&\zr\!\end{bsmallmatrix}\!\}\!\!\!\!\!\!\!\!$ & \multicolumn{4}{c}{$ A_{1}A_{2}A_{1}A_{2}^{2}$} & \\
$\!\!\!\{\!\begin{bsmallmatrix}\!\zr&\mn\!\\\!\pl&\zr\!\end{bsmallmatrix},\!\begin{bsmallmatrix}\!\pl&\pl\!\\\!\mn&\mn\!\end{bsmallmatrix}\!\}\!\!\!\!\!\!\!\!$ & $A_{1}A_{2}\!\!\!$ & &$\!\!\!\{\!\begin{bsmallmatrix}\!\zr&\mn\!\\\!\pl&\pl\!\end{bsmallmatrix},\!\begin{bsmallmatrix}\!\zr&\mn\!\\\!\pl&\mn\!\end{bsmallmatrix}\!\}\!\!\!\!\!\!\!\!$ & $A_{1}A_{2}\!\!\!$ & &$\!\!\!\{\!\begin{bsmallmatrix}\!\zr&\mn\!\\\!\pl&\pl\!\end{bsmallmatrix},\!\begin{bsmallmatrix}\!\pl&\zr\!\\\!\zr&\zr\!\end{bsmallmatrix}\!\}\!\!\!\!\!\!\!\!$ & $A_{1}\!\!\!$ & \!\!(C)\!\!&$\!\!\!\{\!\begin{bsmallmatrix}\!\zr&\mn\!\\\!\pl&\pl\!\end{bsmallmatrix},\!\begin{bsmallmatrix}\!\pl&\zr\!\\\!\zr&\mn\!\end{bsmallmatrix}\!\}\!\!\!\!\!\!\!\!$ & $A_{1}A_{2}\!\!\!$ & &$\!\!\!\{\!\begin{bsmallmatrix}\!\zr&\mn\!\\\!\pl&\pl\!\end{bsmallmatrix},\!\begin{bsmallmatrix}\!\pl&\pl\!\\\!\zr&\zr\!\end{bsmallmatrix}\!\}\!\!\!\!\!\!\!\!$ & $A_{1}\!\!\!$ & \!\!(C)\!\!\\
$\!\!\!\{\!\begin{bsmallmatrix}\!\zr&\mn\!\\\!\pl&\pl\!\end{bsmallmatrix},\!\begin{bsmallmatrix}\!\pl&\pl\!\\\!\zr&\pl\!\end{bsmallmatrix}\!\}\!\!\!\!\!\!\!\!$ & $A_{1}A_{2}^{3}\!\!\!$ & &$\!\!\!\{\!\begin{bsmallmatrix}\!\zr&\mn\!\\\!\pl&\pl\!\end{bsmallmatrix},\!\begin{bsmallmatrix}\!\pl&\pl\!\\\!\zr&\mn\!\end{bsmallmatrix}\!\}\!\!\!\!\!\!\!\!$ & $A_{1}\!\!\!$ & \!\!(C)\!\!&$\!\!\!\{\!\begin{bsmallmatrix}\!\zr&\mn\!\\\!\pl&\pl\!\end{bsmallmatrix},\!\begin{bsmallmatrix}\!\pl&\mn\!\\\!\zr&\zr\!\end{bsmallmatrix}\!\}\!\!\!\!\!\!\!\!$ & $A_{1}^{2}A_{2}\!\!\!$ & &$\!\!\!\{\!\begin{bsmallmatrix}\!\zr&\mn\!\\\!\pl&\pl\!\end{bsmallmatrix},\!\begin{bsmallmatrix}\!\pl&\mn\!\\\!\zr&\pl\!\end{bsmallmatrix}\!\}\!\!\!\!\!\!\!\!$ & $A_{1}^{2}A_{2}^{3}\!\!\!$ & &$\!\!\!\{\!\begin{bsmallmatrix}\!\zr&\mn\!\\\!\pl&\pl\!\end{bsmallmatrix},\!\begin{bsmallmatrix}\!\pl&\mn\!\\\!\zr&\mn\!\end{bsmallmatrix}\!\}\!\!\!\!\!\!\!\!$ & $A_{1}A_{2}\!\!\!$ & \\
$\!\!\!\{\!\begin{bsmallmatrix}\!\zr&\mn\!\\\!\pl&\pl\!\end{bsmallmatrix},\!\begin{bsmallmatrix}\!\pl&\zr\!\\\!\pl&\zr\!\end{bsmallmatrix}\!\}\!\!\!\!\!\!\!\!$ & $A_{1}^{2}A_{2}\!\!\!$ & &$\!\!\!\{\!\begin{bsmallmatrix}\!\zr&\mn\!\\\!\pl&\pl\!\end{bsmallmatrix},\!\begin{bsmallmatrix}\!\pl&\zr\!\\\!\pl&\pl\!\end{bsmallmatrix}\!\}\!\!\!\!\!\!\!\!$ & $A_{1}^{2}A_{2}^{3}\!\!\!$ & &$\!\!\!\{\!\begin{bsmallmatrix}\!\zr&\mn\!\\\!\pl&\pl\!\end{bsmallmatrix},\!\begin{bsmallmatrix}\!\pl&\zr\!\\\!\pl&\mn\!\end{bsmallmatrix}\!\}\!\!\!\!\!\!\!\!$ & $A_{1}A_{2}\!\!\!$ & &$\!\!\!\{\!\begin{bsmallmatrix}\!\zr&\mn\!\\\!\pl&\pl\!\end{bsmallmatrix},\!\begin{bsmallmatrix}\!\pl&\pl\!\\\!\pl&\mn\!\end{bsmallmatrix}\!\}\!\!\!\!\!\!\!\!$ & $A_{2}\!\!\!$ & &\\
$\!\!\!\{\!\begin{bsmallmatrix}\!\zr&\mn\!\\\!\pl&\pl\!\end{bsmallmatrix},\!\begin{bsmallmatrix}\!\pl&\mn\!\\\!\pl&\zr\!\end{bsmallmatrix}\!\}\!\!\!\!\!\!\!\!$ & \multicolumn{4}{c}{$ A_{1}^{2}A_{2}A_{1}A_{2}^{2}A_{1}A_{2}$} & &$\!\!\!\{\!\begin{bsmallmatrix}\!\zr&\mn\!\\\!\pl&\pl\!\end{bsmallmatrix},\!\begin{bsmallmatrix}\!\pl&\mn\!\\\!\pl&\pl\!\end{bsmallmatrix}\!\}\!\!\!\!\!\!\!\!$ & \multicolumn{4}{c}{$ A_{1}A_{2}^{2}A_{1}A_{2}^{3}$} & &$\!\!\!\{\!\begin{bsmallmatrix}\!\zr&\mn\!\\\!\pl&\pl\!\end{bsmallmatrix},\!\begin{bsmallmatrix}\!\pl&\mn\!\\\!\pl&\mn\!\end{bsmallmatrix}\!\}\!\!\!\!\!\!\!\!$ & $A_{1}A_{2}\!\!\!$ & \\
$\!\!\!\{\!\begin{bsmallmatrix}\!\zr&\mn\!\\\!\pl&\pl\!\end{bsmallmatrix},\!\begin{bsmallmatrix}\!\pl&\zr\!\\\!\mn&\zr\!\end{bsmallmatrix}\!\}\!\!\!\!\!\!\!\!$ & $A_{1}\!\!\!$ & \!\!(C)\!\!&$\!\!\!\{\!\begin{bsmallmatrix}\!\zr&\mn\!\\\!\pl&\pl\!\end{bsmallmatrix},\!\begin{bsmallmatrix}\!\pl&\zr\!\\\!\mn&\pl\!\end{bsmallmatrix}\!\}\!\!\!\!\!\!\!\!$ & $A_{1}A_{2}^{3}\!\!\!$ & &$\!\!\!\{\!\begin{bsmallmatrix}\!\zr&\mn\!\\\!\pl&\pl\!\end{bsmallmatrix},\!\begin{bsmallmatrix}\!\pl&\zr\!\\\!\mn&\mn\!\end{bsmallmatrix}\!\}\!\!\!\!\!\!\!\!$ & $A_{1}\!\!\!$ & \!\!(C)\!\!&$\!\!\!\{\!\begin{bsmallmatrix}\!\zr&\mn\!\\\!\pl&\pl\!\end{bsmallmatrix},\!\begin{bsmallmatrix}\!\pl&\pl\!\\\!\mn&\pl\!\end{bsmallmatrix}\!\}\!\!\!\!\!\!\!\!$ & \multicolumn{4}{c}{$ A_{1}A_{2}A_{1}A_{2}A_{1}A_{2}^{2}$} & \\
$\!\!\!\{\!\begin{bsmallmatrix}\!\zr&\mn\!\\\!\pl&\pl\!\end{bsmallmatrix},\!\begin{bsmallmatrix}\!\pl&\pl\!\\\!\mn&\mn\!\end{bsmallmatrix}\!\}\!\!\!\!\!\!\!\!$ & $A_{1}\!\!\!$ & \!\!(C)\!\!&$\!\!\!\{\!\begin{bsmallmatrix}\!\zr&\mn\!\\\!\pl&\pl\!\end{bsmallmatrix},\!\begin{bsmallmatrix}\!\pl&\mn\!\\\!\mn&\mn\!\end{bsmallmatrix}\!\}\!\!\!\!\!\!\!\!$ & $A_{2}\!\!\!$ & &$\!\!\!\{\!\begin{bsmallmatrix}\!\zr&\mn\!\\\!\pl&\mn\!\end{bsmallmatrix},\!\begin{bsmallmatrix}\!\pl&\zr\!\\\!\zr&\zr\!\end{bsmallmatrix}\!\}\!\!\!\!\!\!\!\!$ & $A_{1}\!\!\!$ & \!\!(C)\!\!&$\!\!\!\{\!\begin{bsmallmatrix}\!\zr&\mn\!\\\!\pl&\mn\!\end{bsmallmatrix},\!\begin{bsmallmatrix}\!\pl&\zr\!\\\!\zr&\mn\!\end{bsmallmatrix}\!\}\!\!\!\!\!\!\!\!$ & $A_{1}A_{2}\!\!\!$ & &$\!\!\!\{\!\begin{bsmallmatrix}\!\zr&\mn\!\\\!\pl&\mn\!\end{bsmallmatrix},\!\begin{bsmallmatrix}\!\pl&\pl\!\\\!\zr&\zr\!\end{bsmallmatrix}\!\}\!\!\!\!\!\!\!\!$ & $A_{1}^{2}A_{2}\!\!\!$ & \\
$\!\!\!\{\!\begin{bsmallmatrix}\!\zr&\mn\!\\\!\pl&\mn\!\end{bsmallmatrix},\!\begin{bsmallmatrix}\!\pl&\pl\!\\\!\zr&\pl\!\end{bsmallmatrix}\!\}\!\!\!\!\!\!\!\!$ & $A_{1}^{2}A_{2}^{3}\!\!\!$ & &$\!\!\!\{\!\begin{bsmallmatrix}\!\zr&\mn\!\\\!\pl&\mn\!\end{bsmallmatrix},\!\begin{bsmallmatrix}\!\pl&\pl\!\\\!\zr&\mn\!\end{bsmallmatrix}\!\}\!\!\!\!\!\!\!\!$ & $A_{1}A_{2}\!\!\!$ & &$\!\!\!\{\!\begin{bsmallmatrix}\!\zr&\mn\!\\\!\pl&\mn\!\end{bsmallmatrix},\!\begin{bsmallmatrix}\!\pl&\mn\!\\\!\zr&\zr\!\end{bsmallmatrix}\!\}\!\!\!\!\!\!\!\!$ & $A_{1}\!\!\!$ & \!\!(C)\!\!&$\!\!\!\{\!\begin{bsmallmatrix}\!\zr&\mn\!\\\!\pl&\mn\!\end{bsmallmatrix},\!\begin{bsmallmatrix}\!\pl&\mn\!\\\!\zr&\pl\!\end{bsmallmatrix}\!\}\!\!\!\!\!\!\!\!$ & $A_{1}A_{2}^{3}\!\!\!$ & &$\!\!\!\{\!\begin{bsmallmatrix}\!\zr&\mn\!\\\!\pl&\mn\!\end{bsmallmatrix},\!\begin{bsmallmatrix}\!\pl&\mn\!\\\!\zr&\mn\!\end{bsmallmatrix}\!\}\!\!\!\!\!\!\!\!$ & $A_{1}\!\!\!$ & \!\!(C)\!\!\\
$\!\!\!\{\!\begin{bsmallmatrix}\!\zr&\mn\!\\\!\pl&\mn\!\end{bsmallmatrix},\!\begin{bsmallmatrix}\!\pl&\zr\!\\\!\pl&\zr\!\end{bsmallmatrix}\!\}\!\!\!\!\!\!\!\!$ & $A_{1}\!\!\!$ & \!\!(C)\!\!&$\!\!\!\{\!\begin{bsmallmatrix}\!\zr&\mn\!\\\!\pl&\mn\!\end{bsmallmatrix},\!\begin{bsmallmatrix}\!\pl&\zr\!\\\!\pl&\pl\!\end{bsmallmatrix}\!\}\!\!\!\!\!\!\!\!$ & $A_{1}A_{2}^{3}\!\!\!$ & &$\!\!\!\{\!\begin{bsmallmatrix}\!\zr&\mn\!\\\!\pl&\mn\!\end{bsmallmatrix},\!\begin{bsmallmatrix}\!\pl&\zr\!\\\!\pl&\mn\!\end{bsmallmatrix}\!\}\!\!\!\!\!\!\!\!$ & $A_{1}\!\!\!$ & \!\!(C)\!\!&$\!\!\!\{\!\begin{bsmallmatrix}\!\zr&\mn\!\\\!\pl&\mn\!\end{bsmallmatrix},\!\begin{bsmallmatrix}\!\pl&\pl\!\\\!\pl&\mn\!\end{bsmallmatrix}\!\}\!\!\!\!\!\!\!\!$ & $A_{2}\!\!\!$ & &\\
$\!\!\!\{\!\begin{bsmallmatrix}\!\zr&\mn\!\\\!\pl&\mn\!\end{bsmallmatrix},\!\begin{bsmallmatrix}\!\pl&\mn\!\\\!\pl&\pl\!\end{bsmallmatrix}\!\}\!\!\!\!\!\!\!\!$ & \multicolumn{4}{c}{$ A_{1}A_{2}A_{1}A_{2}A_{1}A_{2}^{2}$} & &$\!\!\!\{\!\begin{bsmallmatrix}\!\zr&\mn\!\\\!\pl&\mn\!\end{bsmallmatrix},\!\begin{bsmallmatrix}\!\pl&\mn\!\\\!\pl&\mn\!\end{bsmallmatrix}\!\}\!\!\!\!\!\!\!\!$ & $A_{1}\!\!\!$ & \!\!(C)\!\!&$\!\!\!\{\!\begin{bsmallmatrix}\!\zr&\mn\!\\\!\pl&\mn\!\end{bsmallmatrix},\!\begin{bsmallmatrix}\!\pl&\zr\!\\\!\mn&\zr\!\end{bsmallmatrix}\!\}\!\!\!\!\!\!\!\!$ & $A_{1}^{2}A_{2}\!\!\!$ & &$\!\!\!\{\!\begin{bsmallmatrix}\!\zr&\mn\!\\\!\pl&\mn\!\end{bsmallmatrix},\!\begin{bsmallmatrix}\!\pl&\zr\!\\\!\mn&\pl\!\end{bsmallmatrix}\!\}\!\!\!\!\!\!\!\!$ & $A_{1}^{2}A_{2}^{3}\!\!\!$ & \\
$\!\!\!\{\!\begin{bsmallmatrix}\!\zr&\mn\!\\\!\pl&\mn\!\end{bsmallmatrix},\!\begin{bsmallmatrix}\!\pl&\zr\!\\\!\mn&\mn\!\end{bsmallmatrix}\!\}\!\!\!\!\!\!\!\!$ & $A_{1}A_{2}\!\!\!$ & &$\!\!\!\{\!\begin{bsmallmatrix}\!\zr&\mn\!\\\!\pl&\mn\!\end{bsmallmatrix},\!\begin{bsmallmatrix}\!\pl&\pl\!\\\!\mn&\zr\!\end{bsmallmatrix}\!\}\!\!\!\!\!\!\!\!$ & \multicolumn{4}{c}{$ A_{1}^{2}A_{2}A_{1}A_{2}^{2}A_{1}A_{2}$} & &$\!\!\!\{\!\begin{bsmallmatrix}\!\zr&\mn\!\\\!\pl&\mn\!\end{bsmallmatrix},\!\begin{bsmallmatrix}\!\pl&\pl\!\\\!\mn&\pl\!\end{bsmallmatrix}\!\}\!\!\!\!\!\!\!\!$ & \multicolumn{4}{c}{$ A_{1}A_{2}^{2}A_{1}A_{2}^{3}$} & \\
$\!\!\!\{\!\begin{bsmallmatrix}\!\zr&\mn\!\\\!\pl&\mn\!\end{bsmallmatrix},\!\begin{bsmallmatrix}\!\pl&\pl\!\\\!\mn&\mn\!\end{bsmallmatrix}\!\}\!\!\!\!\!\!\!\!$ & $A_{1}A_{2}\!\!\!$ & &$\!\!\!\{\!\begin{bsmallmatrix}\!\zr&\mn\!\\\!\pl&\mn\!\end{bsmallmatrix},\!\begin{bsmallmatrix}\!\pl&\mn\!\\\!\mn&\mn\!\end{bsmallmatrix}\!\}\!\!\!\!\!\!\!\!$ & $A_{2}\!\!\!$ & &$\!\!\!\{\!\begin{bsmallmatrix}\!\pl&\zr\!\\\!\zr&\mn\!\end{bsmallmatrix},\!\begin{bsmallmatrix}\!\pl&\mn\!\\\!\pl&\zr\!\end{bsmallmatrix}\!\}\!\!\!\!\!\!\!\!$ & $A_{1}A_{2}\!\!\!$ & &$\!\!\!\{\!\begin{bsmallmatrix}\!\pl&\zr\!\\\!\zr&\mn\!\end{bsmallmatrix},\!\begin{bsmallmatrix}\!\pl&\mn\!\\\!\pl&\mn\!\end{bsmallmatrix}\!\}\!\!\!\!\!\!\!\!$ & $A_{1}A_{2}\!\!\!$ & &$\!\!\!\{\!\begin{bsmallmatrix}\!\pl&\pl\!\\\!\zr&\pl\!\end{bsmallmatrix},\!\begin{bsmallmatrix}\!\pl&\zr\!\\\!\pl&\mn\!\end{bsmallmatrix}\!\}\!\!\!\!\!\!\!\!$ & $A_{1}^{3}A_{2}\!\!\!$ & \\
$\!\!\!\{\!\begin{bsmallmatrix}\!\pl&\pl\!\\\!\zr&\pl\!\end{bsmallmatrix},\!\begin{bsmallmatrix}\!\pl&\pl\!\\\!\pl&\mn\!\end{bsmallmatrix}\!\}\!\!\!\!\!\!\!\!$ & $A_{2}\!\!\!$ & &$\!\!\!\{\!\begin{bsmallmatrix}\!\pl&\pl\!\\\!\zr&\pl\!\end{bsmallmatrix},\!\begin{bsmallmatrix}\!\pl&\mn\!\\\!\pl&\zr\!\end{bsmallmatrix}\!\}\!\!\!\!\!\!\!\!$ & $A_{1}^{3}A_{2}\!\!\!$ & &$\!\!\!\{\!\begin{bsmallmatrix}\!\pl&\pl\!\\\!\zr&\pl\!\end{bsmallmatrix},\!\begin{bsmallmatrix}\!\pl&\mn\!\\\!\pl&\pl\!\end{bsmallmatrix}\!\}\!\!\!\!\!\!\!\!$ & $A_{1}^{2}A_{2}\!\!\!$ & &$\!\!\!\{\!\begin{bsmallmatrix}\!\pl&\pl\!\\\!\zr&\pl\!\end{bsmallmatrix},\!\begin{bsmallmatrix}\!\pl&\mn\!\\\!\pl&\mn\!\end{bsmallmatrix}\!\}\!\!\!\!\!\!\!\!$ & $A_{1}^{4}A_{2}\!\!\!$ & &$\!\!\!\{\!\begin{bsmallmatrix}\!\pl&\pl\!\\\!\zr&\pl\!\end{bsmallmatrix},\!\begin{bsmallmatrix}\!\pl&\zr\!\\\!\mn&\zr\!\end{bsmallmatrix}\!\}\!\!\!\!\!\!\!\!$ & $A_{1}^{5}A_{2}\!\!\!$ & \\
$\!\!\!\{\!\begin{bsmallmatrix}\!\pl&\pl\!\\\!\zr&\pl\!\end{bsmallmatrix},\!\begin{bsmallmatrix}\!\pl&\zr\!\\\!\mn&\pl\!\end{bsmallmatrix}\!\}\!\!\!\!\!\!\!\!$ & $A_{1}^{4}A_{2}^{4}\!\!\!$ & &$\!\!\!\{\!\begin{bsmallmatrix}\!\pl&\pl\!\\\!\zr&\pl\!\end{bsmallmatrix},\!\begin{bsmallmatrix}\!\pl&\zr\!\\\!\mn&\mn\!\end{bsmallmatrix}\!\}\!\!\!\!\!\!\!\!$ & $A_{1}^{3}A_{2}\!\!\!$ & &$\!\!\!\{\!\begin{bsmallmatrix}\!\pl&\pl\!\\\!\zr&\pl\!\end{bsmallmatrix},\!\begin{bsmallmatrix}\!\pl&\pl\!\\\!\mn&\zr\!\end{bsmallmatrix}\!\}\!\!\!\!\!\!\!\!$ & $A_{1}^{3}A_{2}^{2}\!\!\!$ & &$\!\!\!\{\!\begin{bsmallmatrix}\!\pl&\pl\!\\\!\zr&\pl\!\end{bsmallmatrix},\!\begin{bsmallmatrix}\!\pl&\pl\!\\\!\mn&\pl\!\end{bsmallmatrix}\!\}\!\!\!\!\!\!\!\!$ & $A_{1}^{2}A_{2}^{3}\!\!\!$ & &$\!\!\!\{\!\begin{bsmallmatrix}\!\pl&\pl\!\\\!\zr&\pl\!\end{bsmallmatrix},\!\begin{bsmallmatrix}\!\pl&\pl\!\\\!\mn&\mn\!\end{bsmallmatrix}\!\}\!\!\!\!\!\!\!\!$ & $A_{1}^{4}A_{2}\!\!\!$ & \\
$\!\!\!\{\!\begin{bsmallmatrix}\!\pl&\pl\!\\\!\zr&\pl\!\end{bsmallmatrix},\!\begin{bsmallmatrix}\!\pl&\mn\!\\\!\mn&\mn\!\end{bsmallmatrix}\!\}\!\!\!\!\!\!\!\!$ & $A_{2}\!\!\!$ & &$\!\!\!\{\!\begin{bsmallmatrix}\!\pl&\pl\!\\\!\zr&\mn\!\end{bsmallmatrix},\!\begin{bsmallmatrix}\!\pl&\zr\!\\\!\pl&\zr\!\end{bsmallmatrix}\!\}\!\!\!\!\!\!\!\!$ & $A_{1}A_{2}\!\!\!$ & &$\!\!\!\{\!\begin{bsmallmatrix}\!\pl&\pl\!\\\!\zr&\mn\!\end{bsmallmatrix},\!\begin{bsmallmatrix}\!\pl&\zr\!\\\!\pl&\pl\!\end{bsmallmatrix}\!\}\!\!\!\!\!\!\!\!$ & $A_{1}A_{2}^{3}\!\!\!$ & &$\!\!\!\{\!\begin{bsmallmatrix}\!\pl&\pl\!\\\!\zr&\mn\!\end{bsmallmatrix},\!\begin{bsmallmatrix}\!\pl&\zr\!\\\!\pl&\mn\!\end{bsmallmatrix}\!\}\!\!\!\!\!\!\!\!$ & $A_{1}A_{2}\!\!\!$ & &$\!\!\!\{\!\begin{bsmallmatrix}\!\pl&\pl\!\\\!\zr&\mn\!\end{bsmallmatrix},\!\begin{bsmallmatrix}\!\pl&\pl\!\\\!\pl&\mn\!\end{bsmallmatrix}\!\}\!\!\!\!\!\!\!\!$ & $A_{2}\!\!\!$ & \\
$\!\!\!\{\!\begin{bsmallmatrix}\!\pl&\pl\!\\\!\zr&\mn\!\end{bsmallmatrix},\!\begin{bsmallmatrix}\!\pl&\mn\!\\\!\pl&\zr\!\end{bsmallmatrix}\!\}\!\!\!\!\!\!\!\!$ & $A_{1}A_{2}\!\!\!$ & &$\!\!\!\{\!\begin{bsmallmatrix}\!\pl&\pl\!\\\!\zr&\mn\!\end{bsmallmatrix},\!\begin{bsmallmatrix}\!\pl&\mn\!\\\!\pl&\pl\!\end{bsmallmatrix}\!\}\!\!\!\!\!\!\!\!$ & $A_{1}A_{2}^{2}\!\!\!$ & &$\!\!\!\{\!\begin{bsmallmatrix}\!\pl&\pl\!\\\!\zr&\mn\!\end{bsmallmatrix},\!\begin{bsmallmatrix}\!\pl&\mn\!\\\!\pl&\mn\!\end{bsmallmatrix}\!\}\!\!\!\!\!\!\!\!$ & $A_{1}A_{2}\!\!\!$ & &$\!\!\!\{\!\begin{bsmallmatrix}\!\pl&\pl\!\\\!\zr&\mn\!\end{bsmallmatrix},\!\begin{bsmallmatrix}\!\pl&\zr\!\\\!\mn&\zr\!\end{bsmallmatrix}\!\}\!\!\!\!\!\!\!\!$ & $A_{1}\!\!\!$ & &$\!\!\!\{\!\begin{bsmallmatrix}\!\pl&\pl\!\\\!\zr&\mn\!\end{bsmallmatrix},\!\begin{bsmallmatrix}\!\pl&\zr\!\\\!\mn&\pl\!\end{bsmallmatrix}\!\}\!\!\!\!\!\!\!\!$ & $A_{1}A_{2}^{3}\!\!\!$ & \\
$\!\!\!\{\!\begin{bsmallmatrix}\!\pl&\pl\!\\\!\zr&\mn\!\end{bsmallmatrix},\!\begin{bsmallmatrix}\!\pl&\zr\!\\\!\mn&\mn\!\end{bsmallmatrix}\!\}\!\!\!\!\!\!\!\!$ & $A_{1}\!\!\!$ & \!\!(C)\!\!&$\!\!\!\{\!\begin{bsmallmatrix}\!\pl&\pl\!\\\!\zr&\mn\!\end{bsmallmatrix},\!\begin{bsmallmatrix}\!\pl&\pl\!\\\!\mn&\zr\!\end{bsmallmatrix}\!\}\!\!\!\!\!\!\!\!$ & $A_{1}\!\!\!$ & \!\!(C)\!\!&$\!\!\!\{\!\begin{bsmallmatrix}\!\pl&\pl\!\\\!\zr&\mn\!\end{bsmallmatrix},\!\begin{bsmallmatrix}\!\pl&\pl\!\\\!\mn&\pl\!\end{bsmallmatrix}\!\}\!\!\!\!\!\!\!\!$ & $A_{1}A_{2}^{2}\!\!\!$ & &$\!\!\!\{\!\begin{bsmallmatrix}\!\pl&\pl\!\\\!\zr&\mn\!\end{bsmallmatrix},\!\begin{bsmallmatrix}\!\pl&\pl\!\\\!\mn&\mn\!\end{bsmallmatrix}\!\}\!\!\!\!\!\!\!\!$ & $A_{1}\!\!\!$ & &$\!\!\!\{\!\begin{bsmallmatrix}\!\pl&\pl\!\\\!\zr&\mn\!\end{bsmallmatrix},\!\begin{bsmallmatrix}\!\pl&\mn\!\\\!\mn&\mn\!\end{bsmallmatrix}\!\}\!\!\!\!\!\!\!\!$ & $A_{2}\!\!\!$ & \\
$\!\!\!\{\!\begin{bsmallmatrix}\!\pl&\mn\!\\\!\zr&\pl\!\end{bsmallmatrix},\!\begin{bsmallmatrix}\!\pl&\zr\!\\\!\pl&\zr\!\end{bsmallmatrix}\!\}\!\!\!\!\!\!\!\!$ & $A_{1}^{5}A_{2}\!\!\!$ & &$\!\!\!\{\!\begin{bsmallmatrix}\!\pl&\mn\!\\\!\zr&\pl\!\end{bsmallmatrix},\!\begin{bsmallmatrix}\!\pl&\zr\!\\\!\pl&\pl\!\end{bsmallmatrix}\!\}\!\!\!\!\!\!\!\!$ & $A_{1}^{4}A_{2}^{4}\!\!\!$ & &$\!\!\!\{\!\begin{bsmallmatrix}\!\pl&\mn\!\\\!\zr&\pl\!\end{bsmallmatrix},\!\begin{bsmallmatrix}\!\pl&\zr\!\\\!\pl&\mn\!\end{bsmallmatrix}\!\}\!\!\!\!\!\!\!\!$ & $A_{1}^{3}A_{2}\!\!\!$ & &$\!\!\!\{\!\begin{bsmallmatrix}\!\pl&\mn\!\\\!\zr&\pl\!\end{bsmallmatrix},\!\begin{bsmallmatrix}\!\pl&\pl\!\\\!\pl&\mn\!\end{bsmallmatrix}\!\}\!\!\!\!\!\!\!\!$ & $A_{2}\!\!\!$ & &$\!\!\!\{\!\begin{bsmallmatrix}\!\pl&\mn\!\\\!\zr&\pl\!\end{bsmallmatrix},\!\begin{bsmallmatrix}\!\pl&\mn\!\\\!\pl&\zr\!\end{bsmallmatrix}\!\}\!\!\!\!\!\!\!\!$ & $A_{1}^{3}A_{2}^{2}\!\!\!$ & \\
$\!\!\!\{\!\begin{bsmallmatrix}\!\pl&\mn\!\\\!\zr&\pl\!\end{bsmallmatrix},\!\begin{bsmallmatrix}\!\pl&\mn\!\\\!\pl&\pl\!\end{bsmallmatrix}\!\}\!\!\!\!\!\!\!\!$ & $A_{1}^{2}A_{2}^{3}\!\!\!$ & &$\!\!\!\{\!\begin{bsmallmatrix}\!\pl&\mn\!\\\!\zr&\pl\!\end{bsmallmatrix},\!\begin{bsmallmatrix}\!\pl&\mn\!\\\!\pl&\mn\!\end{bsmallmatrix}\!\}\!\!\!\!\!\!\!\!$ & $A_{1}^{4}A_{2}\!\!\!$ & &$\!\!\!\{\!\begin{bsmallmatrix}\!\pl&\mn\!\\\!\zr&\pl\!\end{bsmallmatrix},\!\begin{bsmallmatrix}\!\pl&\zr\!\\\!\mn&\zr\!\end{bsmallmatrix}\!\}\!\!\!\!\!\!\!\!$ & $A_{1}^{2}A_{2}\!\!\!$ & &$\!\!\!\{\!\begin{bsmallmatrix}\!\pl&\mn\!\\\!\zr&\pl\!\end{bsmallmatrix},\!\begin{bsmallmatrix}\!\pl&\zr\!\\\!\mn&\pl\!\end{bsmallmatrix}\!\}\!\!\!\!\!\!\!\!$ & $A_{1}A_{2}\!\!\!$ & &$\!\!\!\{\!\begin{bsmallmatrix}\!\pl&\mn\!\\\!\zr&\pl\!\end{bsmallmatrix},\!\begin{bsmallmatrix}\!\pl&\zr\!\\\!\mn&\mn\!\end{bsmallmatrix}\!\}\!\!\!\!\!\!\!\!$ & $A_{1}^{3}A_{2}\!\!\!$ & \\
$\!\!\!\{\!\begin{bsmallmatrix}\!\pl&\mn\!\\\!\zr&\pl\!\end{bsmallmatrix},\!\begin{bsmallmatrix}\!\pl&\pl\!\\\!\mn&\zr\!\end{bsmallmatrix}\!\}\!\!\!\!\!\!\!\!$ & $A_{1}^{3}A_{2}\!\!\!$ & &$\!\!\!\{\!\begin{bsmallmatrix}\!\pl&\mn\!\\\!\zr&\pl\!\end{bsmallmatrix},\!\begin{bsmallmatrix}\!\pl&\pl\!\\\!\mn&\pl\!\end{bsmallmatrix}\!\}\!\!\!\!\!\!\!\!$ & $A_{1}^{2}A_{2}\!\!\!$ & &$\!\!\!\{\!\begin{bsmallmatrix}\!\pl&\mn\!\\\!\zr&\pl\!\end{bsmallmatrix},\!\begin{bsmallmatrix}\!\pl&\pl\!\\\!\mn&\mn\!\end{bsmallmatrix}\!\}\!\!\!\!\!\!\!\!$ & $A_{1}^{4}A_{2}\!\!\!$ & &$\!\!\!\{\!\begin{bsmallmatrix}\!\pl&\mn\!\\\!\zr&\pl\!\end{bsmallmatrix},\!\begin{bsmallmatrix}\!\pl&\mn\!\\\!\mn&\mn\!\end{bsmallmatrix}\!\}\!\!\!\!\!\!\!\!$ & $A_{2}\!\!\!$ & &$\!\!\!\{\!\begin{bsmallmatrix}\!\pl&\mn\!\\\!\zr&\mn\!\end{bsmallmatrix},\!\begin{bsmallmatrix}\!\pl&\zr\!\\\!\pl&\zr\!\end{bsmallmatrix}\!\}\!\!\!\!\!\!\!\!$ & $A_{1}\!\!\!$ & \\
$\!\!\!\{\!\begin{bsmallmatrix}\!\pl&\mn\!\\\!\zr&\mn\!\end{bsmallmatrix},\!\begin{bsmallmatrix}\!\pl&\zr\!\\\!\pl&\pl\!\end{bsmallmatrix}\!\}\!\!\!\!\!\!\!\!$ & $A_{1}A_{2}^{3}\!\!\!$ & &$\!\!\!\{\!\begin{bsmallmatrix}\!\pl&\mn\!\\\!\zr&\mn\!\end{bsmallmatrix},\!\begin{bsmallmatrix}\!\pl&\zr\!\\\!\pl&\mn\!\end{bsmallmatrix}\!\}\!\!\!\!\!\!\!\!$ & $A_{1}\!\!\!$ & \!\!(C)\!\!&$\!\!\!\{\!\begin{bsmallmatrix}\!\pl&\mn\!\\\!\zr&\mn\!\end{bsmallmatrix},\!\begin{bsmallmatrix}\!\pl&\pl\!\\\!\pl&\mn\!\end{bsmallmatrix}\!\}\!\!\!\!\!\!\!\!$ & $A_{2}\!\!\!$ & &$\!\!\!\{\!\begin{bsmallmatrix}\!\pl&\mn\!\\\!\zr&\mn\!\end{bsmallmatrix},\!\begin{bsmallmatrix}\!\pl&\mn\!\\\!\pl&\zr\!\end{bsmallmatrix}\!\}\!\!\!\!\!\!\!\!$ & $A_{1}\!\!\!$ & \!\!(C)\!\!&$\!\!\!\{\!\begin{bsmallmatrix}\!\pl&\mn\!\\\!\zr&\mn\!\end{bsmallmatrix},\!\begin{bsmallmatrix}\!\pl&\mn\!\\\!\pl&\pl\!\end{bsmallmatrix}\!\}\!\!\!\!\!\!\!\!$ & $A_{1}A_{2}^{2}\!\!\!$ & \\
$\!\!\!\{\!\begin{bsmallmatrix}\!\pl&\mn\!\\\!\zr&\mn\!\end{bsmallmatrix},\!\begin{bsmallmatrix}\!\pl&\mn\!\\\!\pl&\mn\!\end{bsmallmatrix}\!\}\!\!\!\!\!\!\!\!$ & $A_{1}\!\!\!$ & &$\!\!\!\{\!\begin{bsmallmatrix}\!\pl&\mn\!\\\!\zr&\mn\!\end{bsmallmatrix},\!\begin{bsmallmatrix}\!\pl&\zr\!\\\!\mn&\zr\!\end{bsmallmatrix}\!\}\!\!\!\!\!\!\!\!$ & $A_{1}A_{2}\!\!\!$ & &$\!\!\!\{\!\begin{bsmallmatrix}\!\pl&\mn\!\\\!\zr&\mn\!\end{bsmallmatrix},\!\begin{bsmallmatrix}\!\pl&\zr\!\\\!\mn&\pl\!\end{bsmallmatrix}\!\}\!\!\!\!\!\!\!\!$ & $A_{1}A_{2}^{3}\!\!\!$ & &$\!\!\!\{\!\begin{bsmallmatrix}\!\pl&\mn\!\\\!\zr&\mn\!\end{bsmallmatrix},\!\begin{bsmallmatrix}\!\pl&\zr\!\\\!\mn&\mn\!\end{bsmallmatrix}\!\}\!\!\!\!\!\!\!\!$ & $A_{1}A_{2}\!\!\!$ & &$\!\!\!\{\!\begin{bsmallmatrix}\!\pl&\mn\!\\\!\zr&\mn\!\end{bsmallmatrix},\!\begin{bsmallmatrix}\!\pl&\pl\!\\\!\mn&\zr\!\end{bsmallmatrix}\!\}\!\!\!\!\!\!\!\!$ & $A_{1}A_{2}\!\!\!$ & \\
$\!\!\!\{\!\begin{bsmallmatrix}\!\pl&\mn\!\\\!\zr&\mn\!\end{bsmallmatrix},\!\begin{bsmallmatrix}\!\pl&\pl\!\\\!\mn&\pl\!\end{bsmallmatrix}\!\}\!\!\!\!\!\!\!\!$ & $A_{1}A_{2}^{2}\!\!\!$ & &$\!\!\!\{\!\begin{bsmallmatrix}\!\pl&\mn\!\\\!\zr&\mn\!\end{bsmallmatrix},\!\begin{bsmallmatrix}\!\pl&\pl\!\\\!\mn&\mn\!\end{bsmallmatrix}\!\}\!\!\!\!\!\!\!\!$ & $A_{1}A_{2}\!\!\!$ & &$\!\!\!\{\!\begin{bsmallmatrix}\!\pl&\mn\!\\\!\zr&\mn\!\end{bsmallmatrix},\!\begin{bsmallmatrix}\!\pl&\mn\!\\\!\mn&\mn\!\end{bsmallmatrix}\!\}\!\!\!\!\!\!\!\!$ & $A_{2}\!\!\!$ & &$\!\!\!\{\!\begin{bsmallmatrix}\!\pl&\pl\!\\\!\pl&\mn\!\end{bsmallmatrix},\!\begin{bsmallmatrix}\!\pl&\mn\!\\\!\pl&\zr\!\end{bsmallmatrix}\!\}\!\!\!\!\!\!\!\!$ & $A_{1}\!\!\!$ & &$\!\!\!\{\!\begin{bsmallmatrix}\!\pl&\pl\!\\\!\pl&\mn\!\end{bsmallmatrix},\!\begin{bsmallmatrix}\!\pl&\mn\!\\\!\pl&\mn\!\end{bsmallmatrix}\!\}\!\!\!\!\!\!\!\!$ & $A_{1}\!\!\!$ & \\
$\!\!\!\{\!\begin{bsmallmatrix}\!\pl&\mn\!\\\!\pl&\pl\!\end{bsmallmatrix},\!\begin{bsmallmatrix}\!\pl&\mn\!\\\!\pl&\mn\!\end{bsmallmatrix}\!\}\!\!\!\!\!\!\!\!$ & $A_{1}^{2}A_{2}\!\!\!$ & &$\!\!\!\{\!\begin{bsmallmatrix}\!\pl&\mn\!\\\!\pl&\pl\!\end{bsmallmatrix},\!\begin{bsmallmatrix}\!\pl&\zr\!\\\!\mn&\pl\!\end{bsmallmatrix}\!\}\!\!\!\!\!\!\!\!$ & $A_{1}A_{2}^{2}\!\!\!$ & &$\!\!\!\{\!\begin{bsmallmatrix}\!\pl&\mn\!\\\!\pl&\pl\!\end{bsmallmatrix},\!\begin{bsmallmatrix}\!\pl&\zr\!\\\!\mn&\mn\!\end{bsmallmatrix}\!\}\!\!\!\!\!\!\!\!$ & $A_{1}^{2}A_{2}\!\!\!$ & &$\!\!\!\{\!\begin{bsmallmatrix}\!\pl&\mn\!\\\!\pl&\pl\!\end{bsmallmatrix},\!\begin{bsmallmatrix}\!\pl&\pl\!\\\!\mn&\zr\!\end{bsmallmatrix}\!\}\!\!\!\!\!\!\!\!$ & \multicolumn{4}{c}{$ A_{1}^{2}A_{2}A_{1}A_{2}A_{1}A_{2}$} & \\
$\!\!\!\{\!\begin{bsmallmatrix}\!\pl&\mn\!\\\!\pl&\pl\!\end{bsmallmatrix},\!\begin{bsmallmatrix}\!\pl&\pl\!\\\!\mn&\mn\!\end{bsmallmatrix}\!\}\!\!\!\!\!\!\!\!$ & $A_{1}^{2}A_{2}\!\!\!$ & &$\!\!\!\{\!\begin{bsmallmatrix}\!\pl&\mn\!\\\!\pl&\mn\!\end{bsmallmatrix},\!\begin{bsmallmatrix}\!\pl&\zr\!\\\!\mn&\zr\!\end{bsmallmatrix}\!\}\!\!\!\!\!\!\!\!$ & $A_{1}A_{2}\!\!\!$ & &$\!\!\!\{\!\begin{bsmallmatrix}\!\pl&\mn\!\\\!\pl&\mn\!\end{bsmallmatrix},\!\begin{bsmallmatrix}\!\pl&\zr\!\\\!\mn&\pl\!\end{bsmallmatrix}\!\}\!\!\!\!\!\!\!\!$ & $A_{1}A_{2}^{4}\!\!\!$ & &$\!\!\!\{\!\begin{bsmallmatrix}\!\pl&\mn\!\\\!\pl&\mn\!\end{bsmallmatrix},\!\begin{bsmallmatrix}\!\pl&\zr\!\\\!\mn&\mn\!\end{bsmallmatrix}\!\}\!\!\!\!\!\!\!\!$ & $A_{1}A_{2}\!\!\!$ & &$\!\!\!\{\!\begin{bsmallmatrix}\!\pl&\mn\!\\\!\pl&\mn\!\end{bsmallmatrix},\!\begin{bsmallmatrix}\!\pl&\pl\!\\\!\mn&\zr\!\end{bsmallmatrix}\!\}\!\!\!\!\!\!\!\!$ & $A_{1}A_{2}\!\!\!$ & \\
$\!\!\!\{\!\begin{bsmallmatrix}\!\pl&\mn\!\\\!\pl&\mn\!\end{bsmallmatrix},\!\begin{bsmallmatrix}\!\pl&\pl\!\\\!\mn&\pl\!\end{bsmallmatrix}\!\}\!\!\!\!\!\!\!\!$ & $A_{1}A_{2}^{2}\!\!\!$ & &$\!\!\!\{\!\begin{bsmallmatrix}\!\pl&\mn\!\\\!\pl&\mn\!\end{bsmallmatrix},\!\begin{bsmallmatrix}\!\pl&\pl\!\\\!\mn&\mn\!\end{bsmallmatrix}\!\}\!\!\!\!\!\!\!\!$ & $A_{1}A_{2}\!\!\!$ & &$\!\!\!\{\!\begin{bsmallmatrix}\!\pl&\mn\!\\\!\pl&\mn\!\end{bsmallmatrix},\!\begin{bsmallmatrix}\!\pl&\mn\!\\\!\mn&\zr\!\end{bsmallmatrix}\!\}\!\!\!\!\!\!\!\!$ & $A_{2}\!\!\!$ & &$\!\!\!\{\!\begin{bsmallmatrix}\!\pl&\mn\!\\\!\pl&\mn\!\end{bsmallmatrix},\!\begin{bsmallmatrix}\!\pl&\mn\!\\\!\mn&\mn\!\end{bsmallmatrix}\!\}\!\!\!\!\!\!\!\!$ & $A_{2}\!\!\!$ & &
\end{xtabular}
\par
}

\subsection{\texorpdfstring
    {$3\times3$ binary matrices}%
    {3 times 3 binary matrices}%
}
Due to the large number of cases and the need to save space, 
we give the matrices as a pair of two numbers,
where the binary representation of the numbers corresponds to the entries in the matrices.
The rightmost digit in the binary representation corresponds to the entry 
 in the second matrix, third row, third column.
For example, the pair $3/477$, written in the binary system as
$000'000'011_2/111'011'101_2$ corresponds to the matrix pair
{
$
3/477 \triangleq
000'000'011_2/111'011'101_2 \triangleq
\set{
\begin{bsmallmatrix}
0&0&0\\0&0&1\\0&0&1
\end{bsmallmatrix},
\begin{bsmallmatrix}
1&0&1\\1&1&0\\1&1&1
\end{bsmallmatrix}
}.
$
}
Matrix pairs with the same s.m.p.\ are subsumed under one entry. 
For example, the s.m.p.\ $A_2$ of the matrix pair $3/477$ can be found in the section with 
\framebox(30,10){$\mathbf{A_1=3}$}
and the cell with the entry \framebox(60,10){$476-78,\ A_2$}.\\
%
{%
\centering
\tiny%
\begin{xtabular}{|cc|cc|cc|cc|cc|cc|cc|cc|cc|cc|cc|cc|}
\multicolumn{2}{|c|}{ $\mathbf{A_1 = 2}$\Tstrut\Bstrut}
&\!\!77\!\! & $\!\!\!\!\!\!\!A_{2}\!\!\!$ &\!\!92\!\! & $\!\!\!\!\!\!\!A_{1}A_{2}^{5}\!\!\!$ &\!\!93\!\! & $\!\!\!\!\!\!\!A_{2}\!\!\!$ &\!\!101\!\! & $\!\!\!\!\!\!\!A_{2}\!\!\!$ &\!\!104\!\! & $\!\!\!\!\!\!\!A_{1}A_{2}\!\!\!$ &\!\!105\!\! & $\!\!\!\!\!\!\!A_{1}A_{2}^{2}\!\!\!$ &\!\!\!\!108--9\!\!\!\!\!\! & $\!\!\!\!\!\!\!A_{2}\!\!\!$ &\!\!113\!\! & $\!\!\!\!\!\!\!A_{1}A_{2}^{5}\!\!\!$ &\!\!116\!\! & $\!\!\!\!\!\!\!A_{1}A_{2}^{6}\!\!\!$ &\!\!117\!\! & $\!\!\!\!\!\!\!A_{2}\!\!\!$ &\!\!120\!\! & $\!\!\!\!\!\!\!A_{1}A_{2}^{2}\!\!\!$ \\\hline
\!\!121\!\! & $\!\!\!\!\!\!\!A_{1}A_{2}^{3}\!\!\!$ &\!\!124\!\! & $\!\!\!\!\!\!\!A_{1}A_{2}^{4}\!\!\!$ &\!\!125\!\! & $\!\!\!\!\!\!\!A_{1}A_{2}^{7}\!\!\!$ &\!\!141\!\! & $\!\!\!\!\!\!\!A_{2}\!\!\!$ &\!\!\!\!156--57\!\!\!\!\!\! & $\!\!\!\!\!\!\!A_{2}\!\!\!$ &\!\!169\!\! & $\!\!\!\!\!\!\!A_{1}A_{2}^{5}\!\!\!$ &\!\!\!\!172--73\!\!\!\!\!\! & $\!\!\!\!\!\!\!A_{2}\!\!\!$ &\!\!\!\!184--85 & $\!\!\!\!\!\!\!A_{2}\!\!\!$ &\!\!\!\!188--89\!\!\!\!\!\! & $\!\!\!\!\!\!\!A_{2}\!\!\!$ &\!\!\!\!204--5\!\!\!\!\!\! & $\!\!\!\!\!\!\!A_{2}\!\!\!$ &\!\!\!\!220--21 & $\!\!\!\!\!\!\!A_{2}\!\!\!$ &\!\!225\!\! & $\!\!\!\!\!\!\!A_{1}A_{2}^{6}\!\!\!$ \\
\!\!\!\!228--29\!\!\!\!\!\! & $\!\!\!\!\!\!\!A_{2}\!\!\!$ &\!\!232\!\! & $\!\!\!\!\!\!\!A_{2}\!\!\!$ &\!\!233\!\! & $\!\!\!\!\!\!\!A_{1}A_{2}^{4}\!\!\!$ &\!\!\!\!236--37\!\!\!\!\!\! & $\!\!\!\!\!\!\!A_{2}\!\!\!$ &\!\!\!\!240--41\!\!\!\!\!\! & $\!\!\!\!\!\!\!A_{2}\!\!\!$ &\!\!\!\!244--45\!\!\!\!\!\! & $\!\!\!\!\!\!\!A_{2}\!\!\!$ &\!\!248\!\! & $\!\!\!\!\!\!\!A_{2}\!\!\!$ &\!\!249\!\! & $\!\!\!\!\!\!\!A_{1}A_{2}^{7}\!\!\!$ &\!\!\!\!252--53\!\!\!\!\!\! & $\!\!\!\!\!\!\!A_{2}\!\!\!$ &\!\!\!\!332--33\!\!\!\!\!\! & $\!\!\!\!\!\!\!A_{2}\!\!\!$ &\!\!\!\!348--49\!\!\!\!\!\! & $\!\!\!\!\!\!\!A_{2}\!\!\!$ &\!\!352\!\! & $\!\!\!\!\!\!\!A_{2}\!\!\!$ \\
\end{xtabular}%
\par
}

\section*{Acknowledgements}
Thanks to the anonymous referees which pointed out some technical and mathematical mistakes in the paper.

\newcommand{\doi}[1]{\href{https://doi.org/#1}{doi:~#1}}
\newcommand{\arxiv}[1]{\href{https://arxiv.org/abs/#1}{arXiv:~#1}}
\newcommand{\ttilde}{{\raise-1ex\hbox{\textasciitilde}}}


\begin{thebibliography}{99}


\bibitem{BM2002}
Thierry Bousch,Jean Mairesse,
\emph{Asymptotic height optimization for topical IFS, Tetris heaps, and the finiteness conjecture},
J.\ Amer.\ Math.\ Soc., 15 (2002), 77--111,
\doi{10.1090/S0894-0347-01-00378-2}.

\bibitem{BTV2003}
Vincent~D.~Blondel, Jacques~Theys, Alexander~A.~Vladimirov,
2003.
\emph{An elementary counterexample to the finiteness conjecture},
SIAM J.\ Matrix Anal.\ Appl., 24 (2003), 963--970,
\doi{10.1137/S0895479801397846}.

\bibitem{CCGZ2010}
Antonio~Cicone, Nicola~Guglielmi, Stefano~Serra-Capizzano, Marino~Zennaro,
\emph{Finiteness property of pairs of $2\times2$ sign-matrices via real extremal polytope norm},
Linear Alg.\ Appl., 432 (2010) 01, 796--816,
\doi{10.1016/j.laa.2009.09.022}.

\bibitem{DL1992}
Ingrid~Daubechies, Jeffrey\,C.~Lagarias,
\emph{Two-scale difference equations. ii. local regularity, infinite products of matrices and fractals},
SIAM J.\ Math.\ Anal.\ 23 (1992) 4, 1031--1079,
\doi{10.1137/0523059}.

\bibitem{Grip1996}
Gustav Gripenberg,
\emph{Computing the joint spectral radius},
Linear Alg.\ Appl., 234 (1996), 43--60,
\doi{10.1016/0024-3795(94)00082-4}.

\bibitem{Gur95}
Leonid Gurvits,
\emph{Stability of discrete linear inclusion},
Linear Alg.\ Appl.,
231 (1995), 47--85,
\doi{10.1016/0024-3795(95)90006-3}.


\bibitem{GP2013}
Nicola~Guglielmi, Vladimir\,Yu.~Protasov,
\emph{Exact computation of joint spectral characteristics of linear operators},
Found.\ Comput.\ Math., 13 (2013) 1, 37--39,
\doi{10.1007/s10208-012-9121-0}.

\bibitem{GP2016}
Nicola~Guglielmi, Vladimir\,Yu.~Protasov,
\emph{Invariant polytopes of sets of matrices with applications to regularity of wavelets and subdivisions},
SIAM J.\ Matr.\ Anal.\ Appl., 37 (2016) 1,  18--52,
\doi{10.1137/15M1006945}.




\bibitem{HMST2011}
Kevin\,G.~Hare, Ian\,D.~Morris, Nikita~Sidorov, Jacques~Theys,
\emph{An explicit counterexample to the Lagarias--Wang finiteness conjecture},
Adv.\ Math., 226 (2011) 6, 4667-4701,
\doi{10.1016/j.aim.2010.12.012}.

\bibitem{Jung2009}
Raphael~M.~Jungers, 
\emph{The joint spectral radius. Theory and applications}, 
Lecture Notes in Control and Information Sciences (2009), Springer.

\bibitem{JB2008}
Raphael~M.~Jungers,  Vincent~D.~Blondel, 
\emph{On the finiteness property for rational matrices},
Linear Algebra Appl., 428 (2008), 2283--2295,
\doi{10.1016/j.laa.2007.07.007}.



\bibitem{Koz2005}
Victor Kozyakin,
\emph{A dynamical systems construction of a counterexample to the finiteness
conjecture},
Proceedings of the 44th IEEE Conference on Decision and Control, (2005), 2338--2343,
\doi{10.1109/CDC.2005.1582511}.

\bibitem{LW1995}
Jeffrey~C.~Lagarias, Wang Yang,
\emph{The finiteness conjecture for the generalized spectral radius of a set of matrices},
Lin.\ Algebra Appl., 214 (1995), 17--42,
\doi{10.1016/0024-3795(93)00052-2}.

\bibitem{Mej2020}
Thomas~Mejstrik,
\emph{Improved invariant polytope algorithm and applications}, 
 ACM Trans.\ Math.\ Softw., 46 (2020) 3 (29), 1--26,
\doi{10.1145/3408891}.

\bibitem{ttoolboxes}
Thomas~Mejstrik,
\emph{t-toolboxes for Matlab},
Gitlab, (2018),
\href{https://gitlab.com/tommsch/ttoolboxes}%
{gitlab.com/tommsch/ttoolboxes},
2022-02-23.

\bibitem{MOS01}
Bruce E. Moision, Alon Orlitsky,Paul H. Siegel,
\emph{On codes that avoid specified differences},
IEEE Trans.\ Inf.\ Theory, 47 (2001), 433--442.

\bibitem{MP2022}
Thomas~Mejstrik, Vladimir\,Yu.~Protasov,
\emph{Elliptic polytopes and Lyapunov  norms of linear operators},
\arxiv{2107.02610}.



\bibitem{RS1960}
Gian-Carlo~Rota, Gilbert~Strang,
\emph{A note on the joint spectral radius},
Kon.\ Nederl.\ Acad.\ Wet.\ Proc.\ 63 (60).


\bibitem{Wirth2005}
Fabian Wirth,
\emph{The generalized spectral radius is strictly increasing},
Lin.\ Algebra Appl.,  395 (2005), 141--153,
\doi{10.1016/j.laa.2004.07.013}.


\end{thebibliography}
\end{document}